\documentclass[12pt, a4paper]{amsart}
\usepackage{amscd, amssymb}
\usepackage{amsthm}
\def\newspan{\operatorname{span}}

\def\supp{\operatorname{supp}}
\def\ker{\operatorname{ker}}

\def\sup{\operatorname{sup}}

\def\min{\operatorname{min}}

\def\supp{\operatorname{supp}}

\def\interior{\operatorname{int}}
\def\inf{\operatorname{inf}}
\def\s{\operatorname{s}}
\def\Ext{\operatorname{Ext}}
\def\ER{\operatorname{ER}}
\def\BRV{\operatorname{BRV}}
\def\Lar{\operatorname{Lar}}

\def\C{\mathbb{C}}

\def\N{\mathbb{N}}
\def\Z{\mathbb{Z}}
\def\T{\mathbb{T}}

\def\LL{\mathcal{L}}
\def\OO{\mathcal{O}}
\def\KK{\mathcal{K}}

\def\HH{\mathcal{H}}
\def\AA{\mathcal{A}}
\def\XX{\mathcal{X}}

\def\II{\mathcal{I}}

\def\FF{\mathcal{F}}

\def\EE{\mathcal{E}}
\def\JJ{\mathcal{J}}

\def\NO{\mathcal{N}\mathcal{O}}
\def\FE{\mathcal{F}\mathcal{E}}

\newtheorem{thm}{Theorem}[section]
\newtheorem{cor}[thm]{Corollary}

\newtheorem{lemma}[thm]{Lemma}
\newtheorem{prop}[thm]{Proposition}

\theoremstyle{definition}
\newtheorem{definition}[thm]{Definition}

\newtheorem{notation}[thm]{Notation}

\theoremstyle{remark}
\newtheorem{remark}[thm]{Remark}

\newtheorem{acknowledgements}[thm]{Acknowledgements}

\begin{document}

\title[Higher-rank graph algebras and Exel's crossed product]{Realising the $C^*$-algebra of a higher-rank graph as an Exel crossed product}



\author{Nathan Brownlowe}

\address{Nathan Brownlowe, School of Mathematical and Applied Statistics\\
University of Wollongong\\
NSW 2522\\
Australia}
\email{nathanb@uow.edu.au}

\begin{abstract}
We use the boundary-path space of a finitely-aligned $k$-graph $\Lambda$ to construct a compactly-aligned product system $X$, and we show that the graph algebra $C^*(\Lambda)$ is isomorphic to the Cuntz-Nica-Pimsner algebra $\NO(X)$. In this setting, we introduce the notion of a crossed product by a semigroup of partial endomorphisms and partially-defined transfer operators by defining it to be $\NO(X)$. We then compare this crossed product with other definitions in the literature.
\end{abstract}
\thanks{This research has been supported by the Australian Research Council}
\maketitle

\section{Introduction}\label{Introduction}

In \cite{e1}, Exel proposed a new definition for a crossed product of a unital $C^*$-algebra $A$ by an endomorphism $\alpha$. Exel's definition depends not only on $\alpha$, but also on the choice of {\em transfer operator}: a positive continuous linear map $L:A\to A$ satisfying $L(\alpha(a)b)=aL(b)$. We call a triple $(A,\alpha,L)$ an {\em Exel system}. In his motivating example, Exel finds a family of Exel systems whose crossed products model the Cuntz-Krieger algebras \cite{ck}. This marked the first time a crossed product by an endomorphism could successfully model Cuntz-Krieger algebras.

There are two obvious extensions of Exel's construction. Firstly, to a theory of crossed products of {\em non-unital} $C^*$-algebras capable of modeling the directed-graph generalisation of the Cuntz-Krieger algebras \cite{r}. In \cite{brv}, the authors successfully built such a theory, and they realised the graph algebras of locally-finite graphs with no sources as Exel crossed products \cite[Theorem~5.1]{brv}. The crossed product in question was built from the infinite-path space $E^{\infty}$ and the shift map $\sigma$ on $E^{\infty}$. The hypotheses on $E$ ensure that $E^{\infty}$ is locally compact, and $\sigma$ is everywhere defined, and this allows an Exel system to be defined. The other extension of Exel's work is to crossed products by {\em semigroups} of endomorphisms and transfer operators. In \cite{l}, Larsen has a crossed-product construction for dynamical systems $(A,P,\alpha,L)$ in which $P$ is an abelian semigroup, $\alpha$ is an action of $P$ by endomorphisms, and $L$ is an action of $P$ by transfer operators. Exel has also worked in this area with his theory of interaction groups \cite{e2,e3}.

Motivated by these ideas, we construct a semigroup crossed product that can model the $C^*$-algebras of the higher-rank graphs, or $k$-graphs, of Kumjian and Pask \cite{kp}. The only restriction we place on the $k$-graphs $\Lambda$ whose $C^*$-algebras we model is a necessary finitely-aligned hypothesis, so our result applies in the fullest possible generality. This does come at a price, however, as without a locally-finite hypothesis, or a restriction on sources, the space of infinite paths is not locally compact. To get a locally-compact space we need to consider the bigger {\em boundary-path space} $\partial\Lambda$, and on this space the shift maps $\sigma_n$, $n\in\N^k$, will not in general be everywhere defined. This means we can not form Exel systems, or even a dynamical system in the sense of Larsen \cite{l}. We overcome this problem by first ignoring the crossed-product construction, and focusing on building a {\em product system}.

A product system of Hilbert $A$-bimodules over a semigroup $P$ is a semigroup $X=\bigsqcup_{p\in P}X_p$ such that each $X_p$ is a Hilbert $A$-bimodule, and $x\otimes_A y\mapsto xy$ determines an isomorphism of $X_p\otimes_A X_q$ onto $X_{pq}$ for each $p,q\in P$. Fowler introduced such product systems in \cite{f}. Fowler also defined a Cuntz-Pimsner covariance condition for representations of product systems, and introduced the universal $C^*$-algebra $\OO(X)$ for Cuntz-Pimsner covariant representations of $X$. This generalised Pimsner's $C^*$-algebra for a single Hilbert bimodule \cite{p}. In \cite{sy}, Sims and Yeend looked at the problem of associating a $C^*$-algebra to product systems which satisfies a gauge-invariant uniqueness theorem, and noted in particular that Fowler's $\OO(X)$ will not in general do the job. For a large class of semigroups, and a class of product systems called compactly aligned, Sims and Yeend introduced a covariance condition for representations --- called Cuntz-Nica-Pimsner covariance --- and a $C^*$-algebra $\NO(X)$ universal for such representations. A gauge-invariant uniqueness theorem for $\NO(X)$ is proved in \cite{clsv}.

We build from $\partial\Lambda$ and the $\sigma_n$ topological graphs in the sense of Katsura \cite{k}, and then we apply the construction from \cite{k} to get Hilbert $C_0(\partial\Lambda)$-bimodules $X_n$. We glue the bimodules together to form the {\em boundary-path product system} $X$ over $\N^k$. This gives a new class of product systems for which the Cuntz-Nica-Pimsner algebra $\NO(X)$ is tractable. The main result in this paper says that for $\Lambda$ a finitely-aligned $k$-graph, the graph algebra $C^*(\Lambda)$ is isomorphic to $\NO(X)$. A result, we feel, that gives extra credence to Sims and Yeend's construction, at least in the case for the semigroup $\N^k$. We then construct for each $n\in\N^k$ a partial endomorphism $\alpha_n$ on $C_0(\partial\Lambda)$ and a partially-defined transfer operator $L_n$, and we {\em define} the crossed product $C_0(\partial\Lambda)\rtimes_{\alpha,L}\N^k$ to be $\NO(X)$. This gives us our desired result: $C_0(\partial\Lambda)\rtimes_{\alpha,L}\N^k\cong C^*(\Lambda)$.  

We begin with some preliminaries in Section~\ref{preliminaries}. We state some necessary definitions from the $k$-graph literature, and we state the definition of the Cuntz-Krieger algebra of a $k$-graph. We then state the definitions from \cite{sy} needed to make sense of the notion of Cuntz-Nica-Pimsner covariance, and the Cuntz-Nica-Pimsner algebra of a compactly-aligned product system. In Section~\ref{the product system section} we construct from a finitely-aligned $k$-graph $\Lambda$ the boundary-path product system $X$. The proof that $X$ is compactly aligned requires substantial detail, so we leave this result for the appendix. In Section~\ref{the CNP} we prove the existence of a canonical isomorphism $C^*(\Lambda)\to \NO(X)$. In Section~\ref{conn to scp} we introduce the crossed product $C_0(\partial\Lambda)\rtimes_{\alpha,L}\N^k$, and we discuss the relationship between this crossed product and the crossed product in \cite{brv}; Exel and Royer's crossed product by a partial endomorphism \cite{er}; and Larsen's semigroup crossed product \cite{l}.

\section{Preliminaries}\label{preliminaries}

\subsection{$k$-graphs and their Cuntz-Krieger algebras}

A higher-rank graph, or $k$-graph, is a pair $(\Lambda,d)$ consisting of a countable category $\Lambda$ and a degree functor $d:\Lambda\rightarrow\N^k$ satisfying the unique factorisation property: for all $\lambda\in\Lambda$ and $m,n\in\N^k$ with $d(\lambda)=m+n$, there are unique elements $\mu,\nu\in\Lambda$ such that $d(\mu)=m$, $d(\nu)=n$ and $\lambda=\mu\nu$. We now recall some definitions from the $k$-graph literature; for more details see \cite{fmy}. 

For $\lambda,\mu\in\Lambda$ we denote
\[
\Lambda^{\min}(\lambda,\mu):=\{(\alpha,\beta)\in\Lambda\times\Lambda:\lambda\alpha=\mu\beta\text{ and }d(\lambda\alpha)=d(\lambda)\vee d(\mu)\}.
\]
A $k$-graph $\Lambda$ is {\em finitely aligned} if $\Lambda^{\min}(\lambda,\mu)$ is at most finite for all $\lambda,\mu\in\Lambda$. For each $v\in\Lambda^0$ we denote by $v\Lambda:=\{\lambda\in\Lambda:r(\lambda)=v\}$. A subset $E\subseteq v\Lambda$ is {\em exhaustive} if for every $\mu\in v\Lambda$ there exists a $\lambda\in E$ such that $\Lambda^{\min}(\lambda,\mu)\not=\emptyset$. We denote the set of all {\em finite} exhaustive subsets of $\Lambda$ by $\FE(\Lambda)$. We denote by $v\FE(\Lambda)$ the set $\{E\in\FE(\Lambda):E\subseteq v\Lambda\}$. 

For each $m\in(\N\cup\{\infty\})^k$ we get a $k$-graph $\Omega_{k,m}$ through the following construction. The set $\Omega_{k,m}^0:=\{p\in\N^k:p\le m\}$, and
\[
\Omega_k^*:=\{(p,q)\in\Omega_{k,m}^0\times\Omega_{k,m}^0:p\le q\}.
\]
The range map is given by $r(p,q)=p$; the source map by $s(p,q)=q$; and the degree functor by $d(p,q)=q-p$. Composition is given by $(p,q)(q,r)=(p,r)$. 

For $k$-graph $\Lambda$ we define a {\em graph morphism} $x$ to be a degree-preserving functor from $\Omega_{k,m}$ to $\Lambda$. The range and degree maps are extended to all graph morphisms $x:\Omega_{k,m}\to \Lambda$ by setting $r(x):=x(0)$ and $d(x):=m$. We define the {\em boundary-path space} $\partial\Lambda$ to be the set of all graph morphisms $x$ such that for all $n\in\N^k$ with $n\le d(x)$, and for all $E\in x(n)\FE(\Lambda)$, there exists $\lambda\in E$ such that $x(n,n+d(\lambda))=\lambda$. We know from \cite[Lemmas~5.13]{fmy} that if $\lambda\in\Lambda x(0)$, then $\lambda x\in\partial\Lambda$. We know from \cite[Lemma~5.15]{fmy} that for each $v\in\Lambda^0$ there exists $x\in v\partial\Lambda=\{x\in\partial\Lambda:r(x)=v\}$.

We recall from \cite{sy} the following definition.

\begin{definition}\label{CK family}
Let $\Lambda$ be a finitely-aligned $k$-graph. A {\em Cuntz-Krieger family} in a $C^*$-algebra $B$ is a collection $\{t_{\lambda}:\lambda\in\Lambda\}$ of partial isometries in $B$ satisfying
\begin{itemize}
\item[(CK1)] $\{t_v:v\in\Lambda^0\}$ consists of mutually orthogonal projections;  
\item[(CK2)] $t_{\lambda}t_{\mu}=t_{\lambda\mu}$ whenever $s(\lambda)=r(\mu)$;
\item[(CK3)] $t_{\lambda}^*t_{\mu}=\sum_{(\alpha,\beta)\in\Lambda^{\min}(\lambda,\mu)}t_{\alpha}t_{\beta}^*$; and
\item[(CK4)] $\prod_{\lambda\in E}(t_v-t_{\lambda}t_{\lambda}^*)=0$ for every $v\in \Lambda^0$ and $E\in v\FF\EE(\Lambda)$.
\end{itemize} 
\end{definition}

The {\em Cuntz-Krieger algebra}, or {\em graph algebra}, $C^*(\Lambda)$ is the universal $C^*$-algebra generated by a Cuntz-Krieger $\Lambda$-family.


\subsection{Product systems and their Cuntz-Nica-Pimsner algebras}\label{section on ps etc}

In this subsection we state some key definitions from \cite[Sections~2 and 3]{sy}; see \cite{sy} for more details.

Suppose $A$ is a $C^*$-algebra, and $(G,P)$ is a quasi-lattice ordered group in the sense that: $G$ is a discrete group and $P$ is a subsemigroup of $G$; $P\cap P^{-1}=\{e\}$; and with respect to the partial order $p\le q\Longleftrightarrow p^{-1}q\in P$, any two elements $p,q\in G$ which have a common upper bound in $P$ have a least upper bound $p\vee q\in P$. Suppose $X:=\bigcup_{p\in P}X_p$ is a product system of Hilbert $A$-bimodules. For each $p\in P$ and each $x,y\in X_p$ the operator $\Theta_{x,y}:X_p\to X_p$ defined by $\Theta_{x,y}(z):= x\cdot{\langle y,z\rangle}_A$ is adjointable with $\Theta_{x,y}^*=\Theta_{y,x}$. The span $\KK(X_p):=\overline{\newspan}\{\Theta_{x,y}:x,y\in X_p\}$ is a closed two-sided ideal in $\LL(X_p)$ called the {\em algebra of compact operators on} $X_p$. For $p,q\in P$ with $e<p\le q$ there is a homomorphism $\iota_p^q:\LL(X_p)\rightarrow\LL(X_q)$ characterised by
\begin{equation}\label{char of iota map}
\iota_p^q(S)(xy)=(Sx)y\quad\text{for all $x\in X_p,y\in X_{p^{-1}q}$.}
\end{equation}
For $p\not\le q$ we define $\iota_p^q(S)=0_{\LL(X_q)}$ for all $S\in\LL(X_p)$. The product system $X$ is called {\em compactly aligned} if for all $p,q\in P$ such that $p\vee q<\infty$, and for all $S\in\KK(X_p)$ and $T\in\KK(X_q)$, we have $\iota_p^{p\vee q}(S)\iota_q^{p\vee q}(T)\in\KK(X_{p\vee q})$.

A {\em representation} $\psi$ of $X$ in a $C^*$-algebra $B$ is a map $X\to B$ such that
\begin{itemize}
\item[(1)] each $\psi|_{X_p}:=\psi_p:X_p\to B$ is linear, and $\psi_e:A\to B$ is a homomorphism;
\item[(2)] $\psi_p(x)\psi_y(q)=\psi_{pq}(xy)$ for all $p,q\in P$, $x\in X_p$, and $y\in X_q$; and
\item[(3)] $\psi_e({\langle x,y\rangle}_A^p)=\psi_p(x)^*\psi_p(y)$ for all $p\in P$, and $x,y\in X_p$.
\end{itemize}
It follows from Pimsner's results \cite{p} that for each $p\in P$ there is a homomorphism $\psi^{(p)}:\KK(X_p)\to B$ satisfying $\psi^{(p)}(\Theta_{x,y})=\psi_p(x)\psi_p(y)^*$ for all $x,y\in X_p$. A representation $\psi$ of $X$ is {\em Nica-covariant} if for all $p,q\in P$ and all $S\in\KK(X_p), T\in\KK(X_q)$ we have
\[
\psi^{(p)}(S)\psi^{(q)}(T)=
\begin{cases}
\psi^{(p\vee q)}\big(\iota_p^{p\vee q}(S)\iota_q^{p\vee q}(T)\big) & \text{if $p\vee q<\infty$,}\\
0 & \text{otherwise.}
\end{cases}
\]

We denote by $\phi_p$ the homomorphism $A\rightarrow\LL(X_p)$ implementing the left action of $A$ on $X_p$. We define $I_e=A$, and for each $q\in P\setminus\{e\}$ we write $I_q:=\cap_{e<p\le q}\ker\phi_p$. We then denote by $\widetilde{X}_q$ the Hilbert $A$-bimodule
\[
\widetilde{X}_q:=\bigoplus_{p\le q}X_p\cdot I_{p^{-1}q},
\]
and we denote by $\widetilde{\phi}_q$ the homomorphism implementing the left action of $A$ on $\widetilde{X}_q$. The product system $X$ is said to be $\widetilde{\phi}$-{\em injective} if every $\widetilde{\phi}_q$ is injective.

For $p,q\in P$ with $p\not= e$ there is a homomorphism $\widetilde{\iota}_p^q:\LL(X_p)\rightarrow\LL(\widetilde{X}_q)$ determined by $S\mapsto \bigoplus_{r\le q}\iota_p^r(S)$ for all $S\in\LL(X_p)$; and characterised by
\begin{equation}\label{char of iota-tilde map}
(\widetilde{\iota}_p^q(S)x)(r)=\iota_p^r(S)x(r)\quad\text{for all $x\in\widetilde{X}_q$}.
\end{equation} 
A representation $\psi$ of a $\widetilde{\phi}$-injective product system $X$ in a $C^*$-algebra $B$ is {\em Cuntz-Pimsner covariant} if $\sum_{p\in F}\psi^{(p)}(T_p)=0_B$ whenever $F\subset P$ is finite, $T_p\in\KK(X_p)$ for each $p\in F$, and $\sum_{p\in F}\widetilde{\iota}_p^s(T_p)=0$ for large $s$ (see \cite[Definition~3.8]{sy} for the meaning of ``for large $s$''). A representation $\psi$ of a $\tilde{\phi}$-injective product system $X$ is {\em Cuntz-Nica-Pimsner covariant} if it is both Nica  covariant and Cuntz-Pimsner covariant. It is proved in \cite[Proposition~3.12]{sy} that there exists a $C^*$-algebra $\NO(X)$, called the {\em Cuntz-Nica-Pimsner algebra of} $X$, which is universal for Cuntz-Nica-Pimsner covariant representations of $X$. We denote the universal Cuntz-Nica-Pimsner representation by $j_X:X\rightarrow\NO(X)$.


\section{The boundary-path product system of a $k$-graph}\label{the product system section}

Let $\Lambda$ be a finitely-aligned $k$-graph. For $\lambda\in\Lambda$ we denote the set $D_{\lambda}:=\{x\in\partial\Lambda:x(0,d(\lambda))=\lambda\}$. For $n\in\N^k$ we denote
\[
\AA^n:=\{(\lambda,F):\lambda\in\Lambda\text{ with $d(\lambda)\ge n$},\, F\subseteq s(\lambda)\Lambda\text{ a finite set}\},
\]
and $\AA:=\bigcup_{n\in\N^k}\AA^n$. For $(\lambda,F)\in\AA$ we denote $D_{\lambda F}:=\bigcup_{\nu\in F}D_{\lambda\nu}$. It is proved in \cite[Section~5]{fmy} that the family of sets $\{D_{\lambda}\setminus D_{\lambda F}:(\lambda,F)\in\AA\}$ is a basis of compact and open sets for a Hausdorff topology on $\partial\Lambda$, and $\partial\Lambda$ is a locally compact Hausdorff space. For each $n\in\N^k$ we denote $\partial\Lambda^{\ge n}:=\{x\in\partial\Lambda:d(x)\ge n\}$ and $\partial\Lambda^{\not\ge n}:=\partial\Lambda\setminus\partial\Lambda^{\ge n}$. We now use the subsets $\partial\Lambda^{\ge n}$ to construct topological graphs in the sense of Katsura \cite{k,k2}.

\begin{prop}\label{top graph}
Let $n\in\N^k$ with $\partial\Lambda^{\ge n} \not=\emptyset$. Denote by $\sigma_n$ the shift on $\partial\Lambda^{\ge n}$ given by $\sigma_n(x)(m)=x(m+n)$, and $\iota:\partial\Lambda^{\ge n}\rightarrow\partial\Lambda$ the inclusion mapping. Then $E_n:=(\partial\Lambda,\partial\Lambda^{\ge n},\sigma_n,\iota)$ is a topological graph.
\end{prop}

\begin{proof}
We use the definition of convergence given in \cite[Remark~5.6]{fmy}. Let $(x_i)$ be a sequence in $\partial\Lambda^{\not\ge n}$ converging to $x$. If $x\in\partial\Lambda^{\ge n} $, then there exists $j\in\{1,\dots,k\}$ and a subsequence $({x_i}_k)$ of $(x_i)$ such that ${d({x_i}_k)}_j<{d(x)}_j$ for all ${x_i}_k$. This contradicts that $({x_i}_k)$ converges to $x$, so we must have $x\in\partial\Lambda^{\not\ge n} $, and hence $\partial\Lambda^{\not\ge n} $ is closed in $\partial\Lambda$. Hence $\partial\Lambda^{\ge n}$ is locally compact.

Let $x\in\partial\Lambda^{\ge n}$. Then $D_{x(0,n)}$ is an open neighbourhood of $x$, with $D_{x(0,n)}\subseteq\partial\Lambda^{\ge n}$. The map $\sigma_n|_{D_{x(0,n)}}:D_{x(0,n)}\rightarrow D_{s(x(0,n))}$ is a bijection, and $\sigma_n(D_{x(0,n)})=D_{s(x(0,n))}$ is open in $\partial\Lambda$. Now suppose $\lambda\in s(x(0,n))\Lambda$ and $F\subseteq s(\lambda)\Lambda$. Then
\[
\sigma_n|_{D_{x(0,n)}}(D_{x(0,n)\lambda}\setminus D_{x(0,n)\lambda F})= D_{\lambda}\setminus D_{\lambda F}
\]
is open in $D_{s(x(0,n))}$, and
\[
{\left(\sigma_n|_{D_{x(0,n)}}\right)}^{-1}(D_{\lambda}\setminus D_{\lambda F})=D_{x(0,n)\lambda}\setminus D_{x(0,n)\lambda F}
\]
is open in $D_{x(0,n)}$. Hence, $\sigma_n|_{D_{x(0,n)}}$ is continuous and open, and so it is a homeomorphism of $D_{x(0,n)}$ onto $D_{s(x(0,n))}$. Hence $\sigma_n$ is a local homeomorphism. We know that $\iota$ is continuous, so the result follows.
\end{proof}

We now use Katsura's construction \cite{k} to form Hilbert bimodules.  For $f,g\in C_c(\partial\Lambda^{\ge n})$ and $a\in C_0(\partial\Lambda)$, we define
\begin{equation}\label{the right action}
(f\cdot a)(x):= f(x)a(\sigma_n(x))
\end{equation}
and
\begin{equation}\label{the inner product}
{\langle f,g\rangle}_n(x):=\sum_{\sigma_n(y)=x}\overline{f(y)}g(y).
\end{equation}
We complete $C_c(\partial\Lambda^{\ge n})$ under the norm $\|\cdot\|_n$ given by ${\langle\cdot,\cdot\rangle}_n$ to get a Hilbert $C_0(\partial\Lambda)$-module $X_n=X(E_n)$. The formula
\begin{equation}\label{the left action}
(a\cdot f)(x):=a(\iota(x))f(x)=a(x)f(x),
\end{equation}
defines an action of $C_0(\partial\Lambda)$ by adjointable operators on $X_n$, which we denote by $\phi_n:C_0(\partial\Lambda)\rightarrow\LL(X_n)$, and then $X_n$ becomes a Hilbert $C_0(\partial\Lambda)$-bimodule. For $n\in\N^k$ with $\partial\Lambda^{\ge n} =\emptyset$ we set $X_n:=\{0\}$. Note that $X_0=C_0(\partial\Lambda)$.

\begin{prop}\label{the pi map}
Let $m,n\in\N^k$ with $\partial\Lambda^{\ge m} ,\partial\Lambda^{\ge n} \not=\emptyset$. Then the map
\[
\pi:C_c(\partial\Lambda^{\ge m})\times C_c(\partial\Lambda^{\ge n})\rightarrow C_c(\partial\Lambda^{\ge m+n})
\]
given by $\pi(f,g)(x)=f(x)g(\sigma_m(x))$ is a surjective map which induces an isomorphism $\pi_{m,n}:X_m\otimes X_n\to X_{m+n}$ satisfying $\pi_{m,n}(f\otimes g)=f(g\circ\sigma_m)$.
\end{prop}

To prove this proposition we need some results. To state these results we use the following notation.

\begin{notation}
\noindent\textnormal{(a)} Recall from \cite[Definition~3.10]{fmy} that given $\lambda\in\Lambda$ and $E\subseteq r(\lambda)\Lambda$ we denote
\[
\Ext(\lambda;E):=\bigcup_{\nu\in E}\{\alpha\in\Lambda:(\alpha,\beta)\in\Lambda^{\min}(\lambda,\nu)\text{ for some $\beta\in\Lambda$}\}.
\]
For $\lambda,\mu\in\Lambda$ we denote $F(\lambda,\mu):=\Ext(\lambda;\{\mu\})$. Since $\Lambda$ is finitely aligned, $F(\lambda,\mu)$ is a finite subset of $s(\lambda)\Lambda$, and so $(\lambda,F(\lambda,\mu))\in\AA$. We have
\begin{equation}\label{eq for F notation}
D_{\lambda F(\lambda,\mu)}=D_{\mu F(\mu,\lambda)}.
\end{equation}

\smallskip
\noindent\textnormal{(b)} Let $\lambda,\mu\in\Lambda$ and $x\in\partial\Lambda$ with $d(x)\ge d(\lambda)\vee d(\mu)$. Then we denote by $x_{\lambda}^{\mu}$ the path
\[
x_{\lambda}^{\mu}:=x(d(\lambda),d(\lambda)\vee d(\mu)).
\]
\end{notation}

\begin{lemma}\label{set result for being an algebra}
Let $(\lambda,F),(\mu,G)\in\AA$. Then we have
\begin{equation}\label{eq for set result for being an algebra}
(D_{\lambda}\setminus D_{\lambda F})\cap(D_{\mu}\setminus D_{\mu G})
=\bigsqcup_{(\alpha,\beta)\in\Lambda^{\min}(\lambda,\mu)}D_{\lambda\alpha}\setminus D_{\lambda\alpha F_{\alpha}},
\end{equation}
where
\[
F_{\alpha}:=\left(\bigcup_{\nu\in F}F(\lambda\alpha,\lambda\nu)\right)\cup\left(\bigcup_{\xi\in G}F(\lambda\alpha,\mu\xi)\right).
\]
\end{lemma}

\begin{proof}
The factorisation property ensures that the union in (\ref{eq for set result for being an algebra}) is disjoint. 

Let $x\in(D_{\lambda}\setminus D_{\lambda F})\cap(D_{\mu}\setminus D_{\mu G})$. Then $d(x)\ge d(\lambda)\vee d(\mu)$; the pair $(x_{\lambda}^{\mu},x_{\mu}^{\lambda})\in\Lambda^{\min}(\lambda,\mu)$; and $x\in D_{\lambda x_{\lambda}^{\mu}}$. Using (\ref{eq for F notation}) we have
\[
x\in D_{\lambda x_{\lambda}^{\mu} F(\lambda x_{\lambda}^{\mu},\lambda\nu)}=D_{\lambda\nu F(\lambda\nu,\lambda x_{\lambda}^{\mu})}\Longrightarrow x\in D_{\lambda F},
\]
which contradicts $x\in D_{\lambda}\setminus D_{\lambda F}$, so we must have $x\not\in D_{\lambda x_{\lambda}^{\mu}F(\lambda x_{\lambda}^{\mu},\lambda\nu)}$ for all $\nu\in F$. By symmetry, we also have $x\not\in D_{\lambda x_{\lambda}^{\mu}F(\lambda x_{\lambda}^{\mu},\mu\xi)}$ for all $\xi\in G$. Hence $x\in D_{\lambda x_{\lambda}^{\mu}}\setminus D_{\lambda x_{\lambda}^{\mu}F_{ x_{\lambda}^{\mu}}}$.

Now suppose $y$ is an element of the right-hand-side of (\ref{eq for set result for being an algebra}). So there exists $(\alpha,\beta)\in\Lambda^{\min}(\lambda,\mu)$ with $y\in D_{\lambda\alpha}\setminus D_{\lambda\alpha F_{\alpha}}$. We have $y\in D_{\lambda\alpha}\subseteq D_{\lambda}$. Assume $y\in D_{\lambda\nu}$ for some $\nu\in F$. Then $d(y)\ge d(\lambda\alpha)\vee d(\lambda\nu)$; the pair $(y_{\lambda\alpha}^{\lambda\nu},y_{\lambda\nu}^{\lambda\alpha})\in\Lambda^{\min}(\lambda\alpha,\lambda\nu)$; and $y\in D_{\lambda\alpha F(\lambda\alpha,\lambda\nu)}\subseteq D_{\lambda\alpha F_{\alpha}}$. This is a contradiction, and so $y\not\in D_{\lambda\nu}$ for all $\nu\in F$. Hence $y\in D_{\lambda}\setminus D_{\lambda F}$. By symmetry, we also have $y\in D_{\mu}\setminus D_{\mu G}$. Hence $y\in (D_{\lambda}\setminus D_{\lambda F})\cap(D_{\mu}\setminus D_{\mu G})$.
\end{proof}

\begin{lemma}\label{lem about the finite union}
Let $n\in \N^k$ and $(\lambda,F)\in\AA$ with $D_{\lambda}\setminus D_{\lambda F}\subseteq\partial\Lambda^{\ge n}$. Then we have
\begin{equation}\label{the finite, disj union}
D_{\lambda}\setminus D_{\lambda F}=\bigsqcup_{\mu\in s(\lambda)\Lambda^{d(\lambda)\vee n-d(\lambda)}}D_{\lambda\mu}\setminus D_{\lambda\mu \Ext(\mu;F)},
\end{equation}
where $(\lambda\mu,\Ext(\mu;F))\in\AA^n$ for each $\mu\in s(\lambda)\Lambda^{d(\lambda)\vee n-d(\lambda)}$. 
\end{lemma}

\begin{proof}
The factorisation property ensures that the union in (\ref{the finite, disj union}) is disjoint.

Suppose $x\in D_{\lambda}\setminus D_{\lambda F}$, and consider the path $\mu:=x(d(\lambda),d(\lambda)\vee n)\in s(\lambda)\Lambda^{d(\lambda)\vee n-d(\lambda)}$. Then $x\in D_{\lambda\mu}$. If $x\in D_{\lambda\mu \Ext(\mu;F)}$, then there exists $\nu\in F$ and $(\alpha,\beta)\in\Lambda^{\min}(\mu,\nu)$ with $x\in D_{\lambda\mu\alpha}=D_{\lambda\nu\beta}\subseteq D_{\lambda\nu}\subseteq D_{\lambda F}$. But this is a contradiction, and so we must have $x\in D_{\lambda\mu}\setminus D_{\lambda\mu\Ext(\mu;F)}$.

Now, let $y\in D_{\lambda\mu}\setminus D_{\lambda\mu \Ext(\mu;F)}$ for some $\mu\in s(\lambda)\Lambda^{d(\lambda)\vee n-d(\lambda)}$. Then $y\in D_{\lambda}$. If $y\in D_{\lambda\nu}$ for some $\nu\in F$, then the pair
\[
\Big(y_{\lambda\mu}^{\lambda\nu}\big(d(\lambda),d(\lambda)+d(\mu)\vee d(\nu)\big),y_{\lambda\nu}^{\lambda\mu}\big(d(\lambda),d(\lambda)+d(\mu)\vee d(\nu)\big)\Big)\in\Lambda^{\min}(\mu,\nu),
\]
and $y\in D_{\lambda\mu\Ext(\mu;F)}$. This is a contradiction, and so we must have $y\in D_{\lambda}\setminus D_{\lambda F}$.  

Finally, for each $\mu\in s(\lambda)\Lambda^{d(\lambda)\vee n-d(\lambda)}$ the set $\Ext(\mu;F)$ is finite because $F$ is finite and $\Lambda$ is finitely aligned. We obviously have $d(\lambda\mu)\ge n$, and so $(\lambda\mu,\Ext(\mu;F))\in\AA^n$.  
\end{proof}

\begin{proof}[Proof of Proposition~\ref{the pi map}]
To show that $\pi$ is surjective we let $f\in C_c(\partial\Lambda^{\ge m+n})$. For each $x\in\supp f$ there exists $(\lambda,F)\in\AA$ with $x\in D_{\lambda}\setminus D_{\lambda F}\subseteq \partial\Lambda^{\ge m+n}$. So there exists a subset $\JJ\subseteq\AA$ such that $\supp f\subseteq\bigcup_{(\lambda,F)\in\JJ} D_{\lambda}\setminus D_{\lambda F},$ where $D_{\lambda}\setminus D_{\lambda F}\subseteq\partial\Lambda^{\ge m+n}$ for each $(\lambda,F)\in \JJ$. It follows from Lemma~\ref{lem about the finite union} that each $D_{\lambda}\setminus D_{\lambda F}$ is a disjoint union of sets of the form $D_{\mu}\setminus D_{\mu G}$ with $(\mu,G)\in\AA^{m+n}$, and so there exists a subset $\JJ'\subseteq\AA^{m+n}$ such that $\supp f\subseteq\bigcup_{(\mu,G)\in\JJ'} D_{\mu}\setminus D_{\mu G}$, where $D_{\mu}\setminus D_{\mu G}\subseteq\partial\Lambda^{\ge m+n}$ for each $(\mu,G)\in \JJ'$. Since $\supp f$ is compact, there exists a finite number of pairs $(\mu_j,G_j)\in\JJ'$ with $\supp f\subseteq\bigcup_{j=1}^h D_{\mu_j}\setminus D_{\mu_j G_j}$. Now for each $1\le j\le h$ let $\lambda_j:=\mu_j(m,d(\mu_j))$, and consider the function $\XX_{\cup_j D_{\lambda_j}\setminus D_{\lambda_j G_j}}\in C_c(\partial\Lambda^{\ge n})$. Consider also $\tilde{f}\in C_c(\partial\Lambda^{\ge m})$ which is equal to $f$ on $\partial\Lambda^{\ge m+n}$ and zero on the complement. Then we have $\pi(\tilde{f},\XX_{\cup_j D_{\lambda_j}\setminus D_{\lambda_j G_j}})=f$, and so $\pi$ maps onto $C_c(\partial\Lambda^{\ge m+n})$.

Routine calculations show that $\pi$ is bilinear, and so it induces a surjective linear map $\pi_{m,n}:C_c(\partial\Lambda^{\ge m})\odot C_c(\partial\Lambda^{\ge n})\rightarrow C_c(\partial\Lambda^{\ge m+n})$ satisfying $\pi_{m,n}(f\otimes g)(x)=f(x)g(\sigma_m(x))$. It follows immediately from the formulas (\ref{the right action}) and (\ref{the left action}) that $\pi$ preserves the left and right actions.

To see that $\pi_{m,n}$ preserves the inner product, we let $f,h\in C_c(\partial\Lambda^{\ge m})$ and $g,l\in C_c(\partial\Lambda^{\ge n})$. Then for $x\in\partial\Lambda^{\ge m+n}$ we have
\begin{align}
\langle f\otimes g,h\otimes l\rangle(x) &= {\langle {\langle h,f\rangle}_m\cdot g,l\rangle}_n(x)= \sum_{\sigma_n(y)=x}\overline{{\langle h,f\rangle}_m(y)}\overline{g(y)}l(y)\nonumber\\
&= \sum_{\sigma_n(y)=x}\left(\sum_{\sigma_m(z)=y}h(z)\overline{f(z)}\right)\overline{g(y)}l(y)\nonumber\\
&= \sum_{\sigma_{m+n}(z)=x}\overline{g(\sigma_m(z))}l(\sigma_m(z))\overline{f(z)}h(z).\label{eq for zero checking}
\end{align}

Now
\begin{align*}
{\langle\pi_{m,n}(f\otimes g),\pi_{m,n}(h\otimes l)\rangle}_{m+n}(x) &= \sum_{\sigma_{m+n}(z)=x}\overline{\pi_{m,n}(f\otimes g)(z)}\pi_{m,n}(h\otimes l)(z)\\
&= \sum_{\sigma_{m+n}(z)=x}\overline{f(z)}\overline{g(\sigma_m(z))}h(z)l(\sigma_m(z))\\
&= \langle f\otimes g,h\otimes l\rangle(x),
\end{align*}
and so $\pi_{m,n}$ preserves the inner product. Hence it extends to an isomorphism $\pi_{m,n}:X_m\otimes X_n\rightarrow X_{m+n}$.
\end{proof}

\begin{remark}\label{remark for zeros}
Suppose $\partial\Lambda^{\ge m} ,\partial\Lambda^{\ge n} \not= \emptyset$ and $\partial\Lambda^{\ge m+n}=\emptyset$. We claim that $X_m\otimes X_n=\{0\}$. To see this is true, we assume the contrary. Then there exists $f\in C_c(\partial\Lambda^{\ge m})$ and $g\in C_c(\partial\Lambda^{\ge n})$ with $f\otimes g\not=0$. It follows from Equation~(\ref{eq for zero checking}) that 
\[
\langle f\otimes g,f\otimes g\rangle(x) = \sum_{\sigma_{m+n}(z)=x}{|f(z)|}^2{|g(\sigma_m(z))|}^2,
\] 
and this implies 
\begin{align*}
\langle f\otimes g,f\otimes g\rangle\not=0 &\Longleftrightarrow \sigma_{m+n}^{-1}(x)\not=\emptyset\text{ for some $x\in\partial\Lambda$}\\
&\Longleftrightarrow \partial\Lambda^{\ge m+n}\not=\emptyset.
\end{align*}
This is a contradiction, and so we must have $X_m\otimes X_n=\{0\}=X_{m+n}$.

Now suppose that $\partial\Lambda^{\ge m} \not=\emptyset$ and $\partial\Lambda^{\ge n} =\emptyset$. Then we have $\partial\Lambda^{\ge m+n}=\emptyset$, and so $X_n=\{0\}=X_{m+n}$. Then $X_m\otimes X_n=X_m\otimes\{0\}=\{0\}=X_{m+n}$. So we can extend Proposition~\ref{the pi map} to include all $m,n\in\N^k$, and we think of $\pi_{m,n}$ for $m,n$ as in this remark as the trivial map from $\{0\}$ to itself.
\end{remark}

\begin{prop}\label{the product system}
The family $X:=\sqcup_{n\in \N^k}X_n$ of Hilbert bimodules over $C_0(\partial\Lambda)$ with multiplication given by
\begin{equation}\label{def of mult in ps}
xy:=\pi_{m,n}(x\otimes y)
\end{equation}
is a product system over $\N^k$.
\end{prop}

\begin{proof}
We just need to check that $ax=a\cdot x$ and $xa=x\cdot a$ for all $x\in X_n$, $n\in\N^k$ and $a\in C_0(\partial\Lambda)$, but this follows from (\ref{the right action}), (\ref{the left action}) and the definition of multiplication (\ref{def of mult in ps}). 
\end{proof}

We prove that $X$ is compactly aligned in the Appendix.

Given the definition (\ref{def of mult in ps}) of multiplication within $X$, we now have the following restatement of Proposition~\ref{the pi map}. This corollary plays an important role in subsequent sections. 

\begin{cor}\label{imp cor for dense subs}
Let $n\in\N^k$ and $h\in C_c(\partial\Lambda^{\ge n})$. Then for every $l,m\in\N^k$ with $n=l+m$, there exists $f\in C_c(\partial\Lambda^{\ge l})$ and $g\in C_c(\partial\Lambda^{\ge m})$ with $h=fg$.   
\end{cor}


\section{The Cuntz-Nica-Pimsner algebra $\NO(X)$}\label{the CNP}

Recall that we denote by $j_{X}:X\to \NO(X)$ the universal Cuntz-Nica-Pimsner representation of $X$. For each $m\in\N^k$ we denote by $j_{X,m}$ the restriction of $j_{X}$ to $X_m$. For each $\lambda\in\Lambda$ the set $D_\lambda$ is closed and open, and so the characteristic function $\XX_{D_\lambda}\in C_c(\partial\Lambda^{\ge d(\lambda)})\subset X_{d(\lambda)}$.

\begin{thm}\label{the main result}
Let $\Lambda$ be a finitely aligned $k$-graph and $X$ be the associated product system of Hilbert bimodules given in \textup{Proposition}~\textup{\ref{the product system}}. Denote by $\{s_{\lambda}:\lambda\in\Lambda\}$ the universal Cuntz-Krieger $\Lambda$-family in $C^*(\Lambda)$. There exists an isomorphism $\pi: C^*(\Lambda)\rightarrow\NO(X)$ such that $\pi(s_{\lambda})=j_{X,d(\lambda)}(\XX_{D_{\lambda}})$.
\end{thm}

To prove this result we first show that $S:=\{S_{\lambda}:=j_{X,d(\lambda)}(\XX_{D_{\lambda}}):\lambda\in\Lambda\}$ is a set of partial isometries in $\NO(X)$ satisfying (CK1) and (CK2). We use the Nica covariance of $j_{X}$ to show that $S$ satisfies (CK3), and the Cuntz-Pimsner covariance of  $j_{X}$ to show that $S$ satisfies (CK4). The universal property of $C^*(\Lambda)$ then gives us a map $\pi: C^*(\Lambda)\rightarrow\NO(X)$ with $\pi(s_{\lambda})=j_{X}(\XX_{D_{\lambda}})$ for each $\lambda\in\Lambda$. We show that $S$ generates $\NO(X)$, and we use the gauge-invariant uniqueness theorem for $C^*(\Lambda)$ \cite[Theorem~4.2]{rsy} to prove that $\pi$ is injective.

\begin{prop}\label{the first 2 ck relations}
The set $S=\{S_{\lambda}:\lambda\in\Lambda\}$ is a family of partial isometries satisfying \textup{(CK1)} and \textup{(CK2)}.
\end{prop}

\begin{proof}
Let $\lambda\in\Lambda$. Using (\ref{the right action}) and (\ref{the inner product}) we get $\XX_{D_{\lambda}}\cdot{\langle \XX_{D_{\lambda}},\XX_{D_{\lambda}}\rangle}_{d(\lambda)}=\XX_{D_{\lambda}}$, and it follows that $S_{\lambda}S_{\lambda}^*S_{\lambda}=S_{\lambda}$. It follows from the properties of characteristic functions that $\{S_v=j_{X,0}(\XX_{D_v})\}$ is a set of mutually orthogonal projections, thus (CK1) is satisfied. Relation (CK2) follows from the calculation
\begin{align*}
\XX_{D_{\lambda}}\XX_{D_{\mu}}(x) &= \pi_{d(\lambda),d(\mu)}(\XX_{D_{\lambda}}\otimes\XX_{D_{\mu}})(x)\\
&= \XX_{D_{\lambda}}(x)\XX_{D_{\mu}}(\sigma_{d(\lambda)}(x))\\
&= \begin{cases}
1 & \text{if $x(0,d(\lambda))=\lambda$ and $x(d(\lambda),d(\lambda)+d(\mu))=\mu$,}\\
0 & \text{otherwise}
\end{cases}\\
&= \XX_{D_{\lambda\mu}}(x).\qedhere
\end{align*}
\end{proof}

\begin{prop}\label{ck3}
The set $S$ satisfies relation \textup{(CK3)}:
\[
S_{\lambda}^*S_{\mu}=\sum_{(\alpha,\beta)\in\Lambda^{\min}(\lambda,\mu)}S_{\alpha}S_{\beta}^*\quad\text{for all $\lambda,\mu\in\Lambda$.}
\]
\end{prop}

To prove this proposition we need the next result. For $\lambda,\mu\in\Lambda$ with $d(\lambda)=d(\mu)$ we denote by $\Theta_{\lambda,\mu}$ the rank-one operator $\Theta_{\XX_{D_{\lambda}},\XX_{D_{\mu}}}\in\KK(X_{d(\lambda)})$. 

\begin{lemma}\label{the iota maps}
Let $\lambda,\mu\in\Lambda$. Then we have
\begin{equation}\label{eq for iotas}
\iota_{d(\lambda)}^{d(\lambda)\vee d(\mu)}(\Theta_{\lambda,\lambda})\iota_{d(\mu)}^{d(\lambda)\vee d(\mu)}(\Theta_{\mu,\mu}) = \sum_{(\alpha,\beta)\in\Lambda^{\min}(\lambda,\mu)}\Theta_{\lambda\alpha,\mu\beta}.
\end{equation}
\end{lemma}

\begin{proof}
Let $f\in C_c(\partial\Lambda^{\ge d(\mu)})$ and $g\in C_c(\partial\Lambda^{\ge d(\lambda)\vee d(\mu)-d(\mu)})$. We show that the operators in (\ref{eq for iotas}) agree on the product $fg\in C_c(\partial\Lambda^{\ge d(\lambda)\vee d(\mu)})$, and then the result will follow from Corollary~\ref{imp cor for dense subs} and the fact that $C_c(\partial\Lambda^{\ge d(\lambda)\vee d(\mu)})$ is dense in $X_{d(\lambda)\vee d(\mu)}$.  

We know that for each $\mu\in\Lambda$ we have $\Theta_{\mu,\mu}(f)=\XX_{D_{\mu}}\cdot{\langle \XX_{D_{\mu}},f\rangle}_{d(\mu)}$. It follows from a routine calculation using (\ref{the right action}) and (\ref{the inner product}) that $\Theta_{\mu,\mu}(f)=\XX_{D_{\mu}}f$, where $\XX_{D_{\mu}}f$ is a product of functions in $C_c(\partial\Lambda^{\ge d(\mu)})$. It now follows from (\ref{char of iota map}) that
\begin{equation}\label{ex eq}
\iota_{d(\mu)}^{d(\lambda)\vee d(\mu)}(\Theta_{\mu,\mu})(fg)=(\Theta_{\mu,\mu}(f))g=(\XX_{D_{\mu}}f)g.
\end{equation}
We now use Corollary~\ref{imp cor for dense subs} to factor $(\XX_{D_{\mu}}f)g=hl$, where $h\in C_c(\partial\Lambda^{\ge d(\lambda)})$ and $l\in C_c(\partial\Lambda^{\ge d(\lambda)\vee d(\mu)-d(\lambda)})$. For $x\in\partial\Lambda^{\ge d(\lambda)\vee d(\mu)}$ we have
\begin{align*}
\iota_{d(\lambda)}^{d(\lambda)\vee d(\mu)}(\Theta_{\lambda,\lambda})\iota_{d(\mu)}^{d(\lambda)\vee d(\mu)}(\Theta_{\mu,\mu})(fg)(x)&= \iota_{d(\lambda)}^{d(\lambda)\vee d(\mu)}(\Theta_{\lambda,\lambda})(hl)(x)\\
&= (\XX_{D_{\lambda}}h)l(x)\\
&=
\begin{cases}
hl(x) & \text{if $x\in D_{\lambda}$,}\\
0 & \text{otherwise}
\end{cases}\\
&= 
\begin{cases}
fg(x) & \text{if $x\in D_{\lambda}\cap D_{\mu}$,}\\
0 & \text{otherwise.}
\end{cases}
\end{align*}
We know from Lemma~\ref{set result for being an algebra} that $D_{\lambda}\cap D_{\mu}=\sqcup_{(\alpha,\beta)\in\Lambda^{\min}(\lambda,\mu)}D_{\lambda\alpha}$. So we have
\begin{align*}
\iota_{d(\lambda)}^{d(\lambda)\vee d(\mu)}(\Theta_{\lambda,\lambda})\iota_{d(\mu)}^{d(\lambda)\vee d(\mu)}(\Theta_{\mu,\mu})(fg)(x)
&=
\begin{cases}
fg(x) & \text{if $x\in \sqcup_{(\alpha,\beta)\in\Lambda^{\min}(\lambda,\mu)}D_{\lambda\alpha}$,}\\
0 & \text{otherwise.}
\end{cases}\\
&= \left(\sum_{(\alpha,\beta)\in\Lambda^{\min}(\lambda,\mu)}\Theta_{\lambda\alpha,\mu\beta}\right)(fg)(x),
\end{align*}
and the result follows.
\end{proof}

\begin{proof}[Proof of Proposition~\ref{ck3}]
It follows from the Nica covariance of $j_{X}$ that
\begin{align*}
S_{\lambda}S_{\lambda}^*S_{\mu}S_{\mu}^*&= j_{X}^{(d(\lambda))}(\Theta_{\lambda,\lambda})j_{X}^{(d(\mu))}(\Theta_{\mu,\mu})\\
&= j_{X}^{(d(\lambda)\vee d(\mu))}\left(\iota_{d(\lambda)}^{d(\lambda)\vee d(\mu)}(\Theta_{\lambda,\lambda})\iota_{d(\mu)}^{d(\lambda)\vee d(\mu)}(\Theta_{\mu,\mu})\right). 
\end{align*}
It follows from this equation and Lemma~\ref{the iota maps} that 
\begin{align*}
S_{\lambda}\left(\sum_{(\alpha,\beta)\in\Lambda^{\min}(\lambda,\mu)}S_{\alpha}S_{\beta}^*\right)S_\mu^* &= \sum_{(\alpha,\beta)\in\Lambda^{\min}(\lambda,\mu)}S_{\lambda\alpha}{S_{\mu\beta}}^*\\
&= \sum_{(\alpha,\beta)\in\Lambda^{\min}(\lambda,\mu)}j_{X}^{(d(\lambda)\vee d(\mu))}(\Theta_{\lambda\alpha,\mu\beta})\\
&=  j_{X}^{(d(\lambda)\vee d(\mu))}\left(\sum_{(\alpha,\beta)\in\Lambda^{\min}(\lambda,\mu)}\Theta_{\lambda\alpha,\mu\beta}\right)\\
&= S_{\lambda}S_{\lambda}^*S_{\mu}S_{\mu}^*.
\end{align*}
It then follws that
\begin{align*}
S_{\lambda}^*S_{\mu} = (S_{\lambda}^*S_{\lambda}S_{\lambda}^*)(S_{\mu}S_{\mu}^*S_{\mu}) &= S_{\lambda}^*(S_{\lambda}S_{\lambda}^*S_{\mu}S_{\mu}^*)S_{\mu}\\
&= S_{\lambda}^*S_{\lambda}\left(\sum_{(\alpha,\beta)\in\Lambda^{\min}(\lambda,\mu)}S_{\alpha}S_{\beta}^*\right)S_\mu^*S_{\mu}\\
&= \sum_{(\alpha,\beta)\in\Lambda^{\min}(\lambda,\mu)}S_{s(\lambda)}S_{\alpha}{(S_{s(\mu)}S_{\beta})}^*\\
&= \sum_{(\alpha,\beta)\in\Lambda^{\min}(\lambda,\mu)}S_{\alpha}S_{\beta}^*.\qedhere
\end{align*}
\end{proof}

Recall  from Section~\ref{section on ps etc} that $I_n$ is given by $I_n:=\bigcap_{0<m\le n}\ker\phi_m$. To prove that $S$ satisfies (CK4), we need to find families which span dense subspaces of the Hilbert bimodules $X_m\cdot I_{n-m}$, for $m,n\in\N^k$ with $m\le n$. To do this, we must first find families which span dense subspaces of the bimodules $X_n$ and the ideals $I_n$.   

\begin{prop}\label{res for X m}
For each $n\in \N^k$ we have $X_n=\overline{\newspan}\{\XX_{D_{\lambda}\setminus D_{\lambda F}}:(\lambda,F)\in\AA^n\}$.
\end{prop}

\begin{proof}
Let $f\in C_c(\partial\Lambda^{\ge n})$. We can use the same argument as in the beginning of the proof of Proposition~\ref{the pi map} to write $\supp f\subseteq\bigcup_{j=1}^h D_{\mu_j}\setminus D_{\mu_j G_j}$, where $(\mu_j,G_j)\in\AA^n$ and $D_{\mu_j}\setminus D_{\mu_j G_j}\subseteq\partial\Lambda^{\ge n}$ for each $1\le j\le h$. We now take a partition of unity $\rho_1,\dots,\rho_h$ subordinate to $\{D_{\mu_j}\setminus D_{\mu_j G_j}:1\le j\le h\}$, and for $f_j:=f\rho_j\in C(D_{\mu_j}\setminus D_{\mu_j G_j})$ we have
\begin{equation}\label{splitting f}
f=\sum_{j=1}^hf_j.
\end{equation}
Now for each $1\le j\le h$ we have $d(\mu_j)\ge n$. So $\sigma_n$ is injective on $D_{\mu_j}\setminus D_{\mu_j G_j}$, and hence 
\begin{equation}\label{matching norms}
\|f_j\|_n=\sup\{|f_j(x)|:x\in D_{\mu_j}\setminus D_{\mu_j G_j}\}={\|f_j\|}_{\infty}.
\end{equation}
Now, it follows from Lemma~\ref{set result for being an algebra} that for each $(\lambda,F)\in\AA$ the set 
\[
\newspan\{\XX_{D_{\mu}\setminus D_{\mu G}}:(\mu,G)\in\AA\textup{ and }D_{\mu}\setminus D_{\mu G}\subseteq D_{\lambda}\setminus D_{\lambda F}\}
\]
is a subalgbera of $C(D_{\lambda}\setminus D_{\lambda F})$. An application of the Stone-Weierstrass Theorem shows that the closure of that span is equal to $C(D_{\lambda}\setminus D_{\lambda F})$, and hence each $f_j$ can be uniformly approximated by elements in $\newspan\{\XX_{D_{\lambda}\setminus D_{\lambda F}}:d(\lambda)\ge n\}$. It now follows from (\ref{matching norms}) that $f_j$ can be uniformly approximated by elements in $\newspan\{\XX_{D_{\lambda}\setminus D_{\lambda F}}:d(\lambda)\ge n\}$ with respect to $\|\cdot\|_n$, and then (\ref{splitting f}) says that $f$ can be approximated by elements in $\newspan\{\XX_{D_{\lambda}\setminus D_{\lambda F}}:d(\lambda)\ge n\}$ with respect to $\|\cdot\|_n$. The result follows because $C_c(\partial\Lambda\setminus\partial_n)$ is dense in $X_n$ with respect to $\|\cdot\|_n$.
\end{proof}

\begin{definition}\label{condition K i}
Let $i\in\{1,\dots,k\}$ and $e_i$ denote the standard basis element of $\N^k$. We say that $(\lambda,F)\in\AA$ satisfies condition $K(i)$ if
\[
\mu\in s(\lambda)\Lambda\text{ with }d(\mu)\ge e_i\Longrightarrow D_{\mu}\subseteq D_{\nu}\text{ for some $\nu\in F$}.
\]
\end{definition}

\begin{prop}\label{res for I n}
For each $n\in\N^k$ we have
\[
I_n=\overline{\newspan}\{\XX_{D_{\lambda}\setminus D_{\lambda F}}:n_i>0\Longrightarrow{d(\lambda)}_i=0\text{ and $(\lambda,F)$ sat. cond. $K(i)$}\}.
\]
\end{prop}

To prove this proposition we need the following result.

\begin{lemma}\label{kernel for the e i's}
Let $i\in\{1,\dots,k\}$ and $(\lambda,F)\in\AA$. Then $D_{\lambda}\setminus D_{\lambda F}\subseteq\partial\Lambda^{\not\ge e_i}$ if and only if ${d(\lambda)}_i=0$ and $(\lambda,F)$ satisfies condition $K(i)$. Moreover, we have
\begin{equation}\label{eq for ker phi e i}
\ker\phi_{e_i}=\overline{\newspan}\{\XX_{D_{\lambda}\setminus D_{\lambda F}}:(\lambda,F)\in\AA,\,{d(\lambda)}_i=0\text{ and $(\lambda,F)$ satisfies condition $K(i)$}\}.
\end{equation}
\end{lemma}

\begin{proof}
Suppose $D_{\lambda}\setminus D_{\lambda F}\subseteq\partial\Lambda^{\not\ge e_i}$. Then we obviously have ${d(\lambda)}_i=0$. Suppose that $(\lambda,F)$ does not satisfy condition $K(i)$. Then there exists $\mu\in s(\lambda)\Lambda$ with $d(\mu)\ge e_i$, and $x\in D_{\mu}$ with $x\not\in D_{\nu}$ for all $\nu\in F$. Consider the boundary path $\lambda x$. We have ${d(\lambda x)}_i>0$ and $\lambda x\in D_{\lambda}\setminus D_{\lambda F}$. But ${d(\lambda x)}_i>0\Longrightarrow \lambda x \in\partial\Lambda^{\ge e_i}$, and this is a contradiction, so $(\lambda,F)$ satisfies condition $K(i)$.

Now suppose that ${d(\lambda)}_i=0$ and $(\lambda,F)$ satisfies condition $K(i)$. Assume that $D_{\lambda}\setminus D_{\lambda F}\not\subseteq\partial\Lambda^{\not\ge e_i}$, so there exists $x\in D_{\lambda}\setminus D_{\lambda F}$ with $x\in\partial\Lambda^{\ge e_i}$. This implies that ${d(x)}_i>0$. Consider the edge $\mu:=x(d(\lambda),d(\lambda)+e_i)$, which we know exists because ${d(\lambda)}_i=0$. We have $\mu\in s(\lambda)\Lambda\text{ and }d(\mu)=e_i$. The boundary path $\sigma_{d(\lambda)}(x)$ satisfies $\sigma_{d(\lambda)}(x)\in D_{\mu}$ and $\sigma_{d(\lambda)}(x)\not\in D_{\nu}$ for all $\nu\in F$, and so $D_{\mu}\not\subseteq D_{\nu}\text{ for all $\nu\in F$}$. But this contradicts that $(\lambda,F)$ satisfies condition $K(i)$, so we must have $D_{\lambda}\setminus D_{\lambda F}\subseteq\partial\Lambda^{\not\ge e_i}$.

Now, it follows from Lemma~\ref{set result for being an algebra} and an application of the Stone-Weierstrass Theorem for locally compact spaces that for any open subset $U$ of $\partial\Lambda$ we have
\[
C_0(U)=\overline{\newspan}\{\XX_{D_{\lambda}\setminus D_{\lambda F}}:(\lambda,F)\in\AA\text{ and }D_{\lambda}\setminus D_{\lambda F}\subseteq U\}.
\]
It follows that 
\begin{align*}
\ker\phi_{e_i}&=\{a\in C_0(\partial\Lambda):a|_{\partial\Lambda^{\ge e_i}}=0\}\\
&= \{a\in C_0(\partial\Lambda):a|_{\overline{\partial\Lambda^{\ge e_i}}}=0\}\\
&= C_0(\interior\partial\Lambda^{\not\ge e_i})\\
&= \overline{\newspan}\{\XX_{D_{\lambda}\setminus D_{\lambda F}}: (\lambda,F)\in\AA\text{ and }D_{\lambda}\setminus D_{\lambda F}\subseteq\interior\partial\Lambda^{\not\ge e_i} \}\\
&=\overline{\newspan}\{\XX_{D_{\lambda}\setminus D_{\lambda F}}: (\lambda,F)\in\AA\text{ and }D_{\lambda}\setminus D_{\lambda F}\subseteq\partial\Lambda^{\not\ge e_i} \}\\
&= \overline{\newspan}\{\XX_{D_{\lambda}\setminus D_{\lambda F}}:(\lambda,F)\in\AA,\,{d(\lambda)}_i=0,\text{ $(\lambda,F)$ satisfies condition $K(i)$}\}.
\end{align*}\qedhere
\end{proof}

\begin{proof}[Proof of Proposition~\ref{res for I n}]
We have
\begin{align*}
\ker\phi_{n} &= \{a\in C_0(\partial\Lambda):a|_{\partial\Lambda^{\ge n}}=0\}= \{a\in C_0(\partial\Lambda):a|_{\overline{\partial\Lambda^{\ge n}}}=0\}\\
&= C_0(\interior\partial\Lambda^{\ge n}).
\end{align*}
Since $m\le n\Longrightarrow\partial\Lambda^{\not\ge m} \subseteq\partial\Lambda^{\not\ge n} $, it follows that $m\le n\Rightarrow\ker\phi_m\subseteq\ker\phi_n$. Hence $I_n=\bigcap_{\{i:\,n_i>0\}}\ker\phi_{e_i}$, and the result now follows from Lemma~\ref{kernel for the e i's}.
\end{proof}

\begin{notation}\label{note for curly I's}
Let $n\in\N^k$. We define
\begin{align*}
\II(I_n):=\{(\lambda,F)\in\AA:\,& D_{\lambda}\setminus D_{\lambda F}\not=\emptyset\text{ and }\\
& n_i>0\Longrightarrow d(\lambda)_i=0\text{ and $(\lambda,F)$ satisfies condition $K(i)$}\},
\end{align*}
and for $\mu\in\Lambda$ we write $\mu\II(I_n):=\{(\mu\lambda,F):(\lambda,F)\in\II(I_n)\text{ with }s(\mu)=r(\lambda)\}$. The reason for introducing this notation is that we can now write
\[
I_n=\overline{\newspan}\{ \XX_{D_{\lambda}\setminus D_{\lambda F}}:(\lambda,F)\in\II(I_n)\}.
\]
\end{notation}

\begin{prop}\label{X m dot I n-m}
Let $m,n\in\N^k$ with $m\le n$. Then we have
\begin{equation}\label{eq for X m dot I n-m}
X_m\cdot I_{n-m}
=\overline{\newspan}\{\XX_{D_{\lambda}\setminus D_{\lambda F}}: (\lambda,F)\in\AA^m\cap\lambda(0,m)\II(I_{n-m})\}.
\end{equation}
\end{prop}

\begin{proof}
We have $X_m\cdot I_{n-m}= \overline{\newspan}\{x\cdot a:x\in X_m,\,a\in I_{n-m}\}$. To prove that the right-hand side of (\ref{eq for X m dot I n-m}) is contained in the left-hand side, we let $m,n\in\N^k$ with $m\le n$, and suppose $(\lambda,F)\in\AA^m\cap\lambda(0,m)\II(I_{n-m})$. Then $(\lambda(m,d(\lambda)),F)\in \II(I_{n-m})$, and for $x\in\partial\Lambda^{\ge m}$ we have
\begin{align*}
\XX_{D_{\lambda}\setminus D_{\lambda F}}(x) &= \begin{cases}
1 & \text{if $x(0,d(\lambda))=\lambda$ and $x(0,d(\lambda\nu))\not=\lambda\nu$}\\
& \text{for all $\nu\in F$,}\\
0 & \text{otherwise}
\end{cases}\\
&= \XX_{D_{\lambda}}(x)\XX_{D_{\lambda(m,d(\lambda))}\setminus D_{\lambda(m,d(\lambda)) F}}(\sigma_m(x))\\
&= (\XX_{D_{\lambda}}\cdot\XX_{D_{\lambda(m,d(\lambda))}\setminus D_{\lambda(m,d(\lambda)) F}})(x).
\end{align*}
So $\XX_{D_{\lambda}\setminus D_{\lambda F}}=\XX_{D_{\lambda}}\cdot\XX_{D_{\lambda(m,d(\lambda))}\setminus D_{\lambda(m,d(\lambda)) F}}\in X_m\cdot I_{n-m}$, and it follows that 
\[
\overline{\newspan}\{\XX_{D_{\lambda}\setminus D_{\lambda F}}:(\lambda,F)\in\AA^m\cap\lambda(0,m)\II(I_{n-m})\}\subset X_m\cdot I_{n-m}.
\]

It follows from Proposition~\ref{res for I n} and Proposition~\ref{res for X m} that
\[
X_m\cdot I_{n-m}= \overline{\newspan}\{\XX_{D_{\rho}\setminus D_{\rho F}}\cdot\XX_{D_{\tau}\setminus D_{\tau G}}:(\rho,F)\in\AA^m\text{ and }(\tau,G)\in\II(I_{n-m})\}.
\]
So to prove that the left-hand side of (\ref{eq for X m dot I n-m}) is contained in the right-hand side, it suffices to show that for $(\rho,F)\in\AA^m$ and $(\tau,G)\in\II(I_{n-m})$ the product $\XX_{D_{\rho}\setminus D_{\rho F}}\cdot\XX_{D_{\tau}\setminus D_{\tau G}}$ is an element of the right-hand side. Since $\sigma_m^{-1}$ is continuous, the intersection
\begin{equation}\label{int}
(D_{\rho}\setminus D_{\rho F})\cap\sigma_m^{-1}(D_{\tau}\setminus D_{\tau G})
\end{equation}
is an open and compact subset of $D_{\rho}\setminus D_{\rho F}$. Since it is open, we know there exists a subset $\JJ\subseteq\AA^m$ such that $(D_{\rho}\setminus D_{\rho F})\cap\sigma_m^{-1}(D_{\tau}\setminus D_{\tau G})=\bigcup_{(\eta,H)\in\JJ}D_{\eta}\setminus D_{\eta H}$; since it is compact, there is a finite number, say $h$, of pairs $(\eta_j,H_j)\in\JJ$ with
\[
(D_{\rho}\setminus D_{\rho F})\cap\sigma_m^{-1}(D_{\tau}\setminus D_{\tau G})=\bigcup_{j=1}^hD_{\eta_j}\setminus D_{\eta_j H_j},
\] 
We know from Lemma~\ref{set result for being an algebra} that the intersection of sets in the above finite union is a finite, disjoint union of sets of the same form. So it follows that there is a finite number, say $l$, of pairs $(\mu_j,L_j)\in\AA^m$ and constants $c_j$ such that 
\begin{equation}\label{supp of prod}
\XX_{D_{\rho}\setminus D_{\rho F}}\cdot\XX_{D_{\tau}\setminus D_{\tau G}}=\XX_{(D_{\rho}\setminus D_{\rho F})\cap\sigma_m^{-1}(D_{\tau}\setminus D_{\tau G})}=\sum_{j=1}^lc_j\XX_{D_{\mu_j}\setminus D_{\mu_j L_j}}.
\end{equation}
To finish the proof, we need to show that each $(\mu_j,L_j)\in\mu_j(0,m)\II(I_{n-m})$. Suppose $n_i>m_i$ and $d(\mu_j)_i>m_i$. Then for $x\in D_{\mu_j}\setminus D_{\mu_j L_j}$ we have $\sigma_m(x)\in D_{\tau}\setminus D_{\tau G}$ and ${\sigma_m(x)}_i>0$. Since $d(\tau)_i=0$, there exists a path $\alpha:=\sigma_m(x)(d(\tau),d(\tau)+e_i)$ satisfying $\alpha\in s(\tau)\Lambda^{e_i}$. Since $(\tau,G)$ satisfies condition $K(i)$, we have $D_{\alpha}\subseteq D_{\xi}$ for some $\xi\in G$. But this implies that $\sigma_m(x)=\tau\alpha\sigma_m(x)(d(\tau)+e_i,d(x))\in D_{\tau\xi}\subseteq D_{\tau G}$, which contradicts $\sigma_m(x)\in D_{\tau}\setminus D_{\tau G}$. So we must have $d(\mu_j)_i=m_i$.

Now suppose $n_i>m_i$ and there exists an edge $\zeta\in s(\mu_j)\Lambda^{e_i}$ with $D_{\zeta}\not\subseteq D_{\nu}$ for any $\nu\in L_j$. Let $x\in s(\zeta)\partial\Lambda$. Then $\mu_j\zeta x\in D_{\mu_j}\setminus D_{\mu_j L_j}$, which implies  
\begin{equation}\label{imp for lambda j}
\sigma_m(\mu_j\zeta x)\in D_{\tau}\setminus D_{\tau G}.
\end{equation}
Since $d(\tau)_i=0$, there exists a path $\beta:=\sigma_m(\mu_j\zeta x)(d(\tau),d(\tau)+e_i)$ satisfying $\beta\in s(\tau)\Lambda^{e_i}$. Since $(\tau,G)$ satisfies condition $K(i)$, we have $D_{\beta}\subseteq D_{\xi}$ for some $\xi\in G$. But this implies that $\sigma_m(\mu_j\zeta x)=\tau\beta\sigma_m(\mu_j\zeta x)(d(\tau)+e_i,d(x))\in D_{\tau\xi}\subseteq D_{\tau G}$, which contradicts (\ref{imp for lambda j}). So $D_{\zeta}\subseteq D_{\nu}$ for some $\nu\in L_j$, and hence $(\mu_j,L_j)$ satisfies condition $K(i)$.
\end{proof}

\begin{notation}
Let $m,n\in\N^k$ with $m\le n$. We denote
\[
\II(X_m\cdot I_{n-m}) := \{(\lambda,F): D_{\lambda}\setminus D_{\lambda F}\not=\emptyset,\, (\lambda,F)\in\AA^m\cap\lambda(0,m)\II(I_{n-m})\}.
\]
So we have
\[
X_m\cdot I_{n-m}=\overline{\newspan}\{\XX_{D_{\lambda}\setminus D_{\lambda F}}:(\lambda,F)\in\II(X_m\cdot I_{n-m})\}.
\] 
\end{notation}

\begin{prop}\label{ck4}
The set $S=\{S_\lambda:\lambda\in\Lambda\}$ satisfies \textup{(CK4)}:
\[
\prod_{\mu\in\FF}(S_v-S_{\mu}S_{\mu}^*)=0
\]
for all $v\in\Lambda^0$ and all nonempty finite exhaustive sets $\FF\subset r^{-1}(v)$.
\end{prop}

To prove this proposition we need the following results. For a finite subset $G\subset\Lambda$ we denote by $\vee d(G)$ the element $\bigvee_{\mu\in G}d(\mu)$ of $\N^k$.

\begin{lemma}\label{res for ext}
Let $v\in\Lambda^0$ and $\FF\subseteq v\Lambda$ a finite exhaustive set; $n\in\N^k$ with $n\ge\vee d(\FF)$ and $m\in\N^k$ with $m\le n$; and $\lambda\in v\Lambda$ and $F\subseteq s(\lambda)\Lambda$ with $(\lambda,F)\in\II(X_m\cdot I_{n-m})$. Then there exists $\eta\in\FF$ such that $\lambda$ extends $\eta$.
\end{lemma}

\begin{proof}
Suppose $\lambda$ does not extend any element of $\FF$. Since $D_{\lambda}\setminus D_{\lambda F}\not=\emptyset$, there exists a boundary path $x\in D_{\lambda}\setminus D_{\lambda F}$. Since $\FF$ is exhaustive, there exists $\eta\in \FF$ with $x(0,d(\eta))=\eta$. So $x\in D_{\eta}\cap(D_{\lambda}\setminus D_{\lambda F})$, and the pair $(x_{\lambda}^{\eta},x_{\eta}^{\lambda})\in\Lambda^{\min}(\lambda,\eta)$. Since $\lambda$ does not extend $\eta$, there exists $i\in\{1,\dots,k\}$ with ${d(\lambda)}_i<{d(\eta)}_i$, and hence $d(x_{\lambda}^{\eta})\ge e_i$. Since $m_i\le {d(\lambda)}_i<{d(\eta)}_i\le n_i$, we know $(\lambda,F)$ satisfies condition $K(i)$, and hence $D_{x_{\lambda}^{\eta}}\subseteq D_{\nu}$ for some $\nu\in F$. But this implies that $x\in D_{\lambda x_{\lambda}^{\eta}}\subseteq D_{\lambda\nu}$, which contradicts the fact $x\not\in D_{\lambda F}$. So $\lambda$ must extend an element of $\FF$.
\end{proof}

\begin{lemma}\label{res for elt of x tilde}
Suppose $n\in\N^k$ and $\mu\in\Lambda$ with $d(\mu)\le n$. Consider the element $\widetilde{x}$ given by $\widetilde{x}:=(0,\dots,0,\XX_{D_{\lambda}\setminus D_{\lambda F}},0,\dots,0)\in{\widetilde{X}}_n$, where $(\lambda,F)\in\II(X_m\cdot I_{n-m})$ for $m\le n$. Then we have
\[
{\widetilde{\iota}}_{d(\mu)}^n(\Theta_{\mu,\mu})(\widetilde{x})=
\begin{cases}
\widetilde{x} & \text{if $\lambda$ extends $\mu$,}\\
0 & \text{otherwise}.
\end{cases}
\]
\end{lemma}

\begin{proof}
It follows from (\ref{char of iota-tilde map}) that for $r\le n$ we have
\begin{equation}\label{eq 1 for tilde lem}
{\tilde{\iota}}_{d(\mu)}^n(\Theta_{\mu,\mu})(\widetilde{x})(r)= \iota_{d(\mu)}^r(\Theta_{\mu,\mu})(\widetilde{x}(r)) =\begin{cases}
\iota_{d(\mu)}^m(\Theta_{\mu,\mu})(\XX_{D_{\lambda}\setminus D_{\lambda F}}) & \text{if $r=m$,}\\
0 & \text{otherwise.}
\end{cases}
\end{equation}
Now assume $m\ge d(\mu)$. A straightforward calculation shows that 
\begin{equation}\label{eq 2 for tilde lem}
\XX_{D_{\lambda}\setminus D_{\lambda F}}=\XX_{D_{\lambda(d(\mu),d(\lambda))}\setminus D_{\lambda(0,d(\mu))}}\XX_{D_{\lambda(d(\mu),d(\lambda))F}}.
\end{equation}
We also have
\begin{align}
\Theta_{\mu,\mu}(\XX_{D_{\lambda(0,d(\mu))}})(x) &= \left(\XX_{D_{\mu}}\cdot{\langle \XX_{D_{\mu}}^*,\XX_{D_{\lambda(0,d(\mu))}}\rangle}_{d(\mu)}\right)(x)\nonumber\\
&= \XX_{D_{\mu}}(x){\langle \XX_{D_{\mu}}^*,\XX_{D_{\lambda(0,d(\mu))}}\rangle}_{d(\mu)}(\sigma_{d(\mu)}(x))\nonumber\\
&= \begin{cases}
\sum_{\sigma_{d(\mu)}(y)=\sigma_{d(\mu)}(x)}\XX_{D_{\mu}}(y)\XX_{D_{\lambda(0,d(\mu))}}(y) & \text{if $x(0,d(\mu))=\mu$,}\\
0 & \text{otherwise}
\end{cases}\nonumber\\
&= \begin{cases}
1 & \text{if $\lambda(0,d(\mu))=\mu$ and $x(0,d(\mu))=\mu$,}\\
0 & \text{otherwise}
\end{cases}\nonumber\\
&= \begin{cases}
\XX_{D_{\mu}}(x) & \text{if $\lambda$ extends $\mu$,}\\
0 & \text{otherwise}
\end{cases}\nonumber\\
&= \begin{cases}
\XX_{D_{\lambda(0,d(\mu))}}(x) & \text{if $\lambda$ extends $\mu$,}\\
0 & \text{otherwise.}
\end{cases}\label{eq 3 for tilde lem}
\end{align}
It now follows from Equations (\ref{eq 2 for tilde lem}) and (\ref{eq 3 for tilde lem}) that
\begin{align}
\iota_{d(\mu)}^m(\Theta_{\mu,\mu})(\XX_{D_{\lambda}\setminus D_{\lambda F}}) &= \iota_{d(\mu)}^m(\Theta_{\mu,\mu})(\XX_{D_{\lambda(0,d(\mu))}}\XX_{D_{\lambda(d(\mu),d(\lambda))}\setminus D_{\lambda(d(\mu),d(\lambda))F}})\nonumber\\
&= \Theta_{\mu,\mu}(\XX_{D_{\lambda(0,d(\mu))}})\XX_{D_{\lambda(d(\mu),d(\lambda))}\setminus D_{\lambda(d(\mu),d(\lambda))F}}\nonumber\\
&= \begin{cases}
\XX_{D_{\lambda(0,d(\mu))}}\XX_{D_{\lambda(d(\mu),d(\lambda))}\setminus D_{\lambda(d(\mu),d(\lambda))F}} & \text{if $\lambda$ extends $\mu$,}\\
0 & \text{otherwise}
\end{cases}\nonumber\\
&= \begin{cases}
\XX_{D_{\lambda}\setminus D_{\lambda F}} & \text{if $\lambda$ extends $\mu$,}\\
0 & \text{otherwise.}
\end{cases}\label{eq 4 for tilde lem}
\end{align}
Equations (\ref{eq 1 for tilde lem}) and (\ref{eq 4 for tilde lem}) now give the result.
\end{proof}

We are now ready to prove that $S$ satisfies relation (CK4). The proof runs through the main argument from the proof of \cite[Proposition~5.4]{sy}.

\begin{proof}[Proof of Proposition~\ref{ck4}]
Fix $v\in\Lambda^0$ and a finite exhaustive set $\FF\subset v\Lambda$. We must show that
\[
\prod_{\mu\in\FF}(S_v-S_{\mu}S_{\mu}^*)=0.
\]
Recall from \cite{rs} that for a nonempty subset $G$ of $\FF$, $\Lambda^{\min}(G)$ denotes the set $\{\lambda\in\Lambda:d(\lambda)=\vee d(G),\,\lambda\text{ extends $\mu$ for all $\mu\in G$}\}$. Recall also that $\vee\FF:=\bigcup_{G\subset\FF}\Lambda^{\min}(G)$ is finite and is closed under minimal common extensions. We have
\begin{align*}
\prod_{\mu\in\FF}(S_v-S_{\mu}S_{\mu}^*)&=S_v+\sum_{\substack{\emptyset \not= G \subset F \\ \lambda \in \Lambda^{\min}(G)}} {(-1)}^{|G|}S_{\lambda}S_{\lambda}^*\\
&= j_{X}^{(0)}(\Theta_{v,v})+\sum_{\substack{\emptyset \not= G \subset F \\ \lambda \in \Lambda^{\min}(G)}} {(-1)}^{|G|}j_{X}^{(\vee d(G))}(\Theta_{\lambda,\lambda}),
\end{align*}
where the first equation can be obtained through repeated application of (CK3). Since $j_{X}$ is Cuntz-Pimsner covariant, it suffices to show that for each $q\in\N^k$ there exists $r\ge q$ such that for all $s\ge r$, we have
\[
\widetilde{\iota}_0^s(\Theta_{v,v}) + \sum_{\substack{\emptyset \not= G \subset F \\ \lambda \in \Lambda^{\min}(G)}} {(-1)}^{|G|}{\widetilde{\iota}}_{\vee d(G)}^s(\Theta_{\lambda,\lambda})=0.
\]
For this, fix $q\in\N^k$, let $r=q\vee(\vee d(\FF))$ and fix $s\ge r$. It suffices to show that
\begin{equation}\label{eq 1 for ck4}
\Big(\widetilde{\iota}_0^s(\Theta_{v,v}) + \sum_{\substack{\emptyset \not= G \subset F \\ \lambda \in \Lambda^{\min}(G)}} {(-1)}^{|G|}{\widetilde{\iota}}_{\vee d(G)}^s(\Theta_{\lambda,\lambda})\Big)(\widetilde{x})=0,
\end{equation}
where $\widetilde{x}\in {\widetilde{X}}_s$ is given by $\widetilde{x}:=(0,\dots,0,\XX_{D_{\rho}\setminus D_{\rho F}},0,\dots,0)$, for $(\rho,F)\in\II(X_t\cdot I_{s-t})$, $t\le s$. For any $\mu\in\FF$ we have $s\ge d(\mu)$. It then follows from Lemma~\ref{res for elt of x tilde} that
\begin{equation}\label{eq 2 for ck4}
{\widetilde{\iota}}_{d(\mu)}^s(\Theta_{\mu,\mu})(\widetilde{x})=
\begin{cases}
\widetilde{x} & \text{if $\rho$ extends $\mu$,}\\
0 & \text{otherwise}.
\end{cases}
\end{equation}
Fix a nonempty subset $G$ of $\FF$. Then
\[
\Big(\prod_{\mu\in G}{\widetilde{\iota}}_{d(\mu)}^s(\Theta_{\mu,\mu})\Big)(\widetilde{x})=
\begin{cases}
\widetilde{x} & \text{if $\rho$ extends each $\mu\in G$,}\\
0 & \text{otherwise}.
\end{cases}
\]
The factorisation property implies that $\rho$ extends each $\mu\in G$ if and only if there exists $\lambda\in\Lambda^{\min}(G)$ such that $\rho$ extends $\lambda$. The factorisation property also implies that if there does exist such a $\lambda\in\Lambda^{\min}(G)$, then it is necessarily unique. We therefore have
\[
\Big(\prod_{\mu\in G}{\widetilde{\iota}}_{d(\mu)}^s(\Theta_{\mu,\mu})\Big)(\widetilde{x})=\Big(\sum_{\lambda\in\Lambda^{\min}(G)}{\widetilde{\iota}}_{\vee d(G)}^s(\Theta_{\mu,\mu})\Big)(\widetilde{x}).
\]
Since $G$ was an arbitrary subset of $\FF$, we have
\begin{align*}
\Big(\prod_{\mu\in \FF}({\widetilde{\iota}}_{0}^s(\Theta_{v,v})&-{\widetilde{\iota}}_{d(\mu)}^s(\Theta_{\mu,\mu}))\Big)(\widetilde{x})\\
&= \Big({\widetilde{\iota}}_{0}^s(\Theta_{v,v})+\sum_{\emptyset\not= G\subset\FF}\Big({(-1)}^{|G|}\prod_{\mu\in G}{\widetilde{\iota}}_{d(\mu)}^s(\Theta_{\mu,\mu})\Big)\Big)(\widetilde{x})\\
&= \Big(\widetilde{\iota}_0^s(\Theta_{v,v}) + \sum_{\substack{\emptyset \not= G \subset F \\ \lambda \in \Lambda^{\min}(G)}} {(-1)}^{|G|}{\widetilde{\iota}}_{\vee d(G)}^s(\Theta_{\lambda,\lambda})\Big)(\widetilde{x}).
\end{align*}
Now we can apply Lemma~\ref{res for ext} to see that there exists $\eta\in\FF$ such that $\rho$ extends $\eta$. It now follows from Equation~(\ref{eq 2 for ck4}) that
\begin{align*}
&\Big(\prod_{\mu\in \FF}({\widetilde{\iota}}_{0}^s(\Theta_{v,v})-{\widetilde{\iota}}_{d(\mu)}^s(\Theta_{\mu,\mu}))\Big)(\widetilde{x})\\
&\qquad\qquad= \Big(\prod_{\mu\in \FF\setminus\{\eta\}}({\widetilde{\iota}}_{0}^s(\Theta_{v,v})-{\widetilde{\iota}}_{d(\mu)}^s(\Theta_{\mu,\mu}))\Big)(({\widetilde{\iota}}_{0}^s(\Theta_{v,v})-{\widetilde{\iota}}_{d(\eta)}^s(\Theta_{\eta,\eta})))(\widetilde{x})\\
&\qquad\qquad=0,
\end{align*}
and hence Equation~(\ref{eq 1 for ck4}) is established.
\end{proof}

\begin{proof}[Proof of Theorem~\ref{the main result}]
Lemma~\ref{the first 2 ck relations}, Proposition~\ref{ck3} and Proposition~\ref{ck4} show that the set $S:=\{S_{\lambda}=j_{X}(\XX_{D_{\lambda}}):\lambda\in\Lambda\}$ is a family of partial isometries satisfying the Cuntz-Krieger relations (CK1)--(CK4). It follows from the universal property of $C^*(\Lambda)$ that there exists a homomorphism $\pi:C^*(\Lambda)\rightarrow\NO(X)$ such that $\pi(s_{\lambda})=j_{X}(\XX_{D_{\lambda}})$ for each $\lambda\in\Lambda$. We know from \cite[Proposition~3.12]{sy} that $\NO(X)=\overline{\newspan}\{j_{X}(x)j_{X}(y)^*:x,y\in X\}$. For each $\lambda\in\Lambda$ and $F\subseteq s(\lambda)\Lambda$ we have  $\XX_{D_{\lambda}\setminus D_{\lambda F}}=\XX_{D_{\lambda}}-\sum_{\nu\in F}\XX_{D_{\lambda\nu}}$, and so
\[
j_{X}(\XX_{D_{\lambda}\setminus D_{\lambda F}})=j_{X}(\XX_{D_{\lambda}})-j_{X}\left(\sum_{\nu\in F}\XX_{D_{\lambda\nu}}\right)=S_{\lambda}-\sum_{\nu\in F}S_{\lambda\nu}.
\]
It then follows from Proposition~\ref{res for X m} that $S$ generates $\NO(X)$, and hence $\pi$ is surjective. It follows from \cite[Lemma~5.13(2) and Lemma~5.15]{fmy} that each $D_{\lambda}\not=\emptyset$, and hence each $\XX_{D_{\lambda}}\not=0$. It then follows from \cite[Theorem~4.1]{sy} that each $S_{\lambda}\not=0$. (Note that the quasi-lattice ordered group $(\N^k,\Z^k)$ satisfies \cite[Condition~(3.5)]{sy}, and so \cite[Theorem~4.1]{sy} can indeed be applied.) Since $\pi$ intertwines the gauge actions of $\T^k$ on $\NO(X)$ and $C^*(\Lambda)$, the gauge-invariant uniqueness theorem for $C^*(\Lambda)$ \cite[Theorem~4.2]{rsy} implies that $\pi$ is an isomorphism.
\end{proof}


\section{Connections to semigroup crossed products}\label{conn to scp} 

We begin this section by building a crossed product from a finitely-aligned $k$-graph $\Lambda$. For each $n\in\N^k$ we define a {\em partial endomorphism} $\alpha_n:C_0(\partial\Lambda)\to C_0(\partial\Lambda^{\ge n})$ given by $\alpha_n(f)=f\circ\sigma_n$. We claim that for $f\in C_c(\partial\Lambda^{\ge n})$ the function $L_n(f)$ given by 
\[
L_n(f)(x)=
\begin{cases}
\sum_{\sigma_n(y)=x}f(y) & \text{if $x\in\sigma_n(\partial\Lambda)$,}\\
0 & \text{otherwise}\\
\end{cases}
\]
is well-defined and is an element of $C_c(\partial\Lambda)$. We can cover $\supp f$ with finitely many sets $U_i$ such that $\sigma_n(U_i)$ is open, $\overline{\sigma_n(U_i)}$ is compact, and $\sigma_n|_{U_i}$ is a homeomorphism. The function $f$ must be zero on all but a finite number of points in $\sigma_n^{-1}(x)$. Then near any $x\in\sigma_n(\partial\Lambda)$, $L_n(f)=\sum_{\{i:x\in\sigma_n(U_i)\}} f\circ{(\sigma_n|_{U_i})}^{-1}$ is a finite sum of continuous functions with compact support. Since $\sigma_n(x)$ is open, $L_n(f)\in C_c(\partial\Lambda)$, and the claim is proved. Routine calculations show that each $L_n$ satisfies the transfer-operator identity: $L_n(\alpha_n(f)g)=fL_n(g)$ for all $f\in C_0(\partial\Lambda)$, $g\in C_c(\partial\Lambda^{\ge n})$. Adapting Exel's construction of a Hilbert bimodule \cite{e1} to accommodate the partial maps, and applying it to $(C_0(\partial\Lambda),\alpha_n,L_n)$, gives the Hilbert $C_0(\partial\Lambda)$-bimodule $X_n$ from Section~\ref{the product system section}. So we consider the boundary-path product system $X$, and take the suggested route of \cite[Section~9]{brv} for defining a crossed product for the system $(C_0(\partial\Lambda),\N^k,\alpha,L)$:     

\begin{definition}\label{def of cross prod}
Let $\Lambda$ be a finitely-aligned $k$-graph, and consider the product system $X$ given in Proposition~\ref{the product system}. We define the {\em crossed product} $C_0(\partial\Lambda)\rtimes_{\alpha,L}\N^k$ to be the Cuntz-Nica-Pimsner algebra $\NO(X)$.
\end{definition}

\begin{cor}\label{goal result}
Let $\Lambda$ be a finitely-aligned $k$-graph. Then $C_0(\partial\Lambda)\rtimes_{\alpha,L}\N^k\cong C^*(\Lambda)$.
\end{cor} 

For the remainder of this section we discuss the relationship between the crossed product $C_0(\partial\Lambda)\rtimes_{\alpha,L}\N^k$ and the other crossed products in the literature which are given via transfer operators; namely, the non-unital version of Exel's crossed product \cite{brv}, Exel and Royer's crossed product by a partial endomorphism \cite{er}, and Larsen's crossed product for semigroups \cite{l}. The upshot of this discussion is that, when these crossed products can be defined, they coincide with $C_0(\partial\Lambda)\rtimes_{\alpha,L}\N^k$. To be make things clear, we use the following notation.

\begin{notation}
(1) For $(A,\beta,\LL)$ a dynamical system in the sense of Exel and Royer \cite{er} we denote by $A\rtimes_{\beta,\LL}^{\ER}\N$ the crossed product given in \cite[Definition~1.6]{er}.

\smallskip

(2) For $(A,\beta,\LL)$ a dynamical system in the sense of \cite{brv,e1} we denote by $A\rtimes_{\beta,\LL}^{\BRV}\N$ the crossed product given in \cite[Section~4]{brv}.

\smallskip

(3) For $P$ an abelian semigroup and $(A,P,\beta,\LL)$ a dynamical system in the sense of Larsen \cite{l} we denote by $A\rtimes_{\beta,\LL}^{\Lar}P$ the crossed product given in \cite[Definition~2.2]{l}.
\end{notation}

\subsection{Directed graphs}

Suppose $\Lambda$ is a 1-graph. Then for each $\lambda,\mu\in\Lambda$ we have $|\Lambda^{\min}(\lambda,\mu)|\in\{0,1\}$, and so $\Lambda$ is finitely aligned. As shown in \cite[Examples~10.1--10.2]{r}, $\Lambda$ is the path category of the directed graph $E:=(\Lambda^0,d^{-1}(1),r,s)$. We know from \cite[Proposition~B.1]{rsy} that $C^*(\Lambda)$ coincides with the graph algebra $C^*(E)$ as given in \cite{flr}. We denote by $E^*$ the set of finite paths in $E$ and by $E^{\infty}$ the set of infinite paths in $E$. We define $E_{\inf}^*:=\{\mu\in E^*:|r^{-1}(s(\mu))|=\infty\}$ and $E_{\s}^*:=\{\mu\in E^*:r^{-1}(s(\mu))=\emptyset\},$ so $E_{\inf}^*$ is the set of paths whose source is an infinite receiver, and $E_{\s}^*$ is the set of paths whose source is a source in $E$. Then the boundary-path space $\partial\Lambda$ coincides with $\partial E:=E^{\infty}\cup E_{\inf}^*\cup E_{\s}^*$. We now freely use directed graphs $E$ in place of 1-graphs $\Lambda$ in Definition~\ref{def of cross prod}.

\begin{prop}\label{rel to ER for 1-graph}
Let $E$ be a directed graph. Then $(C_0(\partial E),\alpha, L)$ is a dynamical system in the sense of \cite{er}, and we have $C_0(\partial E)\rtimes_{\alpha,L}\N\cong C_0(\partial E)\rtimes_{\alpha,L}^{\ER}\N$.
\end{prop}   

To prove this proposition we need the following result.

\begin{prop}\label{ER cros prod is kat O_M}
Let $(A,\beta,\LL)$ be a dynamical system in the sense of \cite{er}, and consider the Hilbert $A$-bimodule $M$ constructed in \cite[Section~1]{er}. Then $A\rtimes_{\beta,\LL}^{\ER}\N$ is isomorphic to Katsura's Cuntz-Pimsner algebra $\OO_M$ \cite{kat}.
\end{prop}

\begin{proof}
The arguments in \cite[Section~3]{br} (or \cite[Section~4]{brv}) extend across to this setting, except $A\rtimes_{\beta,\LL}^{\ER}\N$ is defined by modding out redundancies $(a,k)$ with $a\in (\ker\phi)^{\perp}\cap\phi^{-1}(\KK(M))$ instead of $\overline{A\alpha(A)A}\cap\phi^{-1}(\KK(M))$. But $(\ker\phi)^{\perp}\cap\phi^{-1}(\KK(M))$ is precisely the ideal involved in Katsura's definition of $\OO_M$ \cite[Definition~3.5]{kat}.
\end{proof}

\begin{proof}[Proof of Proposition~\ref{rel to ER for 1-graph}]
The construction of the Hilbert $A$-bimodule $M$ from \cite{er} gives $X_1$. We know from \cite[Proposition~5.3]{sy} that $\NO(X)$ is isomorphic to Katsura's $\OO_{X_1}$. We know from Proposition~\ref{ER cros prod is kat O_M} that $C_0(\partial E)\rtimes_{\alpha,L}^{\ER}\N\cong\OO_M$. So we have
\[
C_0(\partial E)\rtimes_{\alpha,L}\N=\NO(X)\cong\OO_{X_1}=\OO_M\cong C_0(\partial E)\rtimes_{\alpha,L}^{\ER}\N.\qedhere
\]
\end{proof}

\subsection{Locally-finite directed graphs with no sources}

For a locally-finite directed graph $\Lambda:=E$ with no sources we have $\partial E=E^{\infty}$. We denote by $\sigma$ the backward shift on $E^{\infty}$, and $\alpha_E$ the endomorphism of $C_0(E^{\infty})$ given by $\alpha_E(f)=f\circ\sigma$. So $\alpha_E=\alpha_1$. For each $f\in C_0(E^{\infty})$ we denote by $L_E(f)$ the function given by
\[
L_E(f)(x)=
\begin{cases}
\frac{1}{|\sigma^{-1}(x)|}\sum_{\sigma(y)=x}f(y) & \text{if $x\in\sigma(E^{\infty})$,}\\
0 & \text{otherwise}\\
\end{cases}
\]
So $L_E$ is the normalised version of $L_1$. It is proved in \cite[Section~2.1]{brv} that $L_E$ is a transfer operator for $(C_0(E^{\infty}),\alpha_E)$. 

\begin{prop}\label{rel to BRV for 1-graph}
Let $E$ be a locally-finite directed graph with no sources. Then we have $C_0(E^{\infty})\rtimes_{\alpha,L}\N\cong C_0(E^{\infty})\rtimes_{\alpha_E,L_E}^{\BRV}\N$.
\end{prop}

\begin{proof}
Recall the construction of the Hilbert $C_0(E^{\infty})$-bimodule $M_{L_E}$ \cite[Section~3]{brv}, and in particular that $q:C_0(E^{\infty})\to M_{L_E}$ denotes the quotient map. Since $E$ is locally finite, the shift $\sigma$ is proper. We can use this fact to find for each $x\in E^{\infty}$ an open neighbourhood $V$ of $\sigma(x)$ such that $|\sigma^{-1}(v)|=|\sigma^{-1}(\sigma(x))|$ for each $v\in V$, and it follows that the map $d:E^{\infty}\to\C$ given by $d(x)=\sqrt{|\sigma^{-1}(\sigma(x))|}$ is continuous. Straightforward calculations show that $U:C_c(E^{\infty})\to M_{L_E}$ given by $U(f)=q(df)$ extends to an isomorphism of $X_1$ onto $M_{L_E}$. So $\OO_{X_1}\cong\OO_{M_{L_E}}$. Since $E$ has no sources, the homomorphism $\phi:C_0(E^{\infty})\to\LL(M_{L_E})$ giving the left action on $M_{L_E}$ is injective, and so ${(\ker\phi)}^{\perp}=C_0(E^{\infty})$. It then follows from \cite[Corollary~4.2]{brv} that $C_0(E^{\infty})\rtimes_{\alpha_E,L_E}^{\BRV}\N\cong\OO_{M_{L_E}}$. Finally, we know from \cite[Proposition~5.3]{sy} that $\NO(X)\cong \OO_{X_1}$, so we have
\[
C_0(E^{\infty})\rtimes_{\alpha,L}\N=\NO(X)\cong\OO_{X_1}\cong\OO_{M_{L_E}}\cong C_0(E^{\infty})\rtimes_{\alpha_E,L_E}^{\BRV}\N.\qedhere
\]
\end{proof}

\subsection{Regular $k$-graphs}
We now examine how $C_0(\partial\Lambda)\rtimes_{\alpha,L}\N^k$ fits in with the theory of Larsen's semigroup crossed products \cite{l}.

If $\Lambda$ is a row-finite $k$-graph with no sources, then $\partial\Lambda$ is the set $\Lambda^{\infty}$ of all graph morphisms from $\Omega_{k,(\infty,\dots,\infty)}$ to $\Lambda$, and the shift maps are everywhere defined. So $\alpha$ is an action by endomorphisms. We say a $k$-graph $\Lambda$ is {\em regular} if it is row-finite with no sources, and there exists $M_1,\dots,M_k\in\N\setminus\{0\}$ such that for each $i\in\{1,\dots,k\}$ we have $|\Lambda^{e_i}v|=M_i$ for all $v\in\Lambda^0$. For each $x\in\Lambda^{\infty}$ and $n\in\N^k$ define
\[
\omega(n,x):={|\sigma_n^{-1}(\sigma_n(x))|}^{-1}=\prod_{i=1}^kM_i^{-n_i}.
\]
Then for each $f\in C_0(\Lambda^{\infty})$ the map $\LL_n(f)$ given by
\[
\LL_{n}(f)(x)= 
\begin{cases}
\sum_{\sigma_{n}(y)=x}\omega(n,y)f(y) & \text{if $x\in\sigma_{n}(\Lambda^{\infty})$,}\\
0 & \text{otherwise}\\
\end{cases}
\]
is a transfer operator for $(C_0(\Lambda^{\infty}),\alpha_n)$. Simple calculations show that 
\[
\sum_{\sigma_n(y)=x}\omega(n,y)=1
\]
for all $x\in\Lambda^{\infty},n\in\N^k$, and that $\omega(m+n,x)=\omega(m,x)\omega(n,\sigma_m(x))$ for all $x\in\Lambda^{\infty},m,n\in\N^k$. Hence \cite[Proposition~2.2]{e-Ren}, which still holds in the non-unital setting, gives an action $\LL$ of $\N^k$ of transfer operators on $C_0(\Lambda^{\infty})$. It follows that $(C_0(\Lambda^{\infty}),\N^k,\alpha,\LL)$ is a dynamical system in the sense of Larsen \cite[Section~2]{l}.

\begin{prop}\label{rel to Lar for k-graph}
Let $\Lambda$ be a regular $k$-graph. Then we have $C_0(\Lambda^{\infty})\rtimes_{\alpha,L}\N^k\cong C_0(\Lambda^{\infty})\rtimes_{\alpha,\LL}^{\Lar}\N^k$.
\end{prop}     

\begin{proof}
We apply the construction in \cite[Section~3.2]{l} to the dynamical system $(C_0(\Lambda^{\infty}),\N^k,\alpha,\LL)$ to form a product system $M=\cup_{n\in\N^k}M_{\LL_n}$, and then \cite[Proposition~4.3]{l} says $C_0(\Lambda^{\infty})\rtimes_{\alpha,\LL}^{\Lar}\N^k$ is isomorphic to Fowler's Cuntz-Pimsner algebra $\OO(M)$ \cite[Proposition~2.9]{f}. Suppose $M_1,\dots,M_k\in\N\setminus\{0\}$ such that for each $i\in\{1,\dots,k\}$ we have $|\Lambda^{e_i}v|=M_i$ for all $v\in\Lambda^0$. For each $n\in\N^k$ denote $M_n:=\prod_{i=1}^kM_i^{-n_i}$. Then the map $f\mapsto q_n(\sqrt{M_n}f)$ from $C_c(\Lambda^{\infty})$ to $M_{\LL_n}$ extends to an isomorphism of $X_n$ onto $M_{\LL_n}$. These maps induce an isomorphism of the product systems $X$ and $M$ (observe the formulae for multiplication within $X$, Proposition~\ref{the pi map}, and $M$, \cite[Equation~3.8]{l}). So $\OO(X)\cong\OO(M)$. 

Recall that each $X_n$ is constructed from the topological graph $(\Lambda^{\infty},\Lambda^{\infty},\sigma_n,\iota)$, where $\iota$ is the inclusion map. It then follows from \cite[Proposition~1.24]{k} that each $\phi_n$ is injective and acts by compact operators. So we can apply \cite[Corollary~5.2]{sy} to see that $\NO(X)$ coincides with $\OO(X)$. So we have
\[
C_0(\Lambda^{\infty})\rtimes_{\alpha,L}\N^k=\NO(X)=\OO(X)\cong\OO(M)\cong C_0(\Lambda^{\infty})\rtimes_{\alpha,\LL}^{\Lar}\N^k.\qedhere
\]
\end{proof}

\subsection{Conclusion} The results in this section justify our decision to define the crossed product $C_0(\partial\Lambda)\rtimes_{\alpha,L}\N^k$ to be the Cuntz-Nica-Pimsner algebra $\NO(X)$, and we propose that the same definition is made for a general crossed product by a quasi-lattice ordered semigroup of partial endomorphisms and partially-defined transfer operators. The problem is that Sims and Yeend's Cuntz-Nica-Pimsner algebra is only appropriate for a particular family (containing $\N^k$) of quasi-lattice ordered semigroups. The ``correct'' definition of a Cuntz-Pimsner algebra of a product system over an arbitrary quasi-lattice ordered semigroup is yet to be found. (See \cite{sy,clsv} for more discussion.)


\section{Appendix}

Recall that for $(G,P)$ a quasi-lattice ordered group, and $X$ a product system over $P$ of Hilbert bimodules, we say that $X$ is {\em compactly aligned} if for all $p,q\in P$ such that $p\vee q<\infty$, and for all $S\in\KK(X_p)$ and $T\in\KK(X_q)$, we have $\iota_p^{p\vee q}(S)\iota_q^{p\vee q}(T)\in\KK(X_{p\vee q})$.

\begin{prop}\label{comp al}
The product system $X$ constructed in \textup{Section}~\textup{\ref{the product system section}} is compactly aligned.
\end{prop}

We start with a definition and some notation.

\begin{definition}\label{disjoint set}
Let $n\in\N^k$. We say that a subset $\JJ\subseteq\AA^n$ is {\em disjoint} if 
\[
(\lambda,F),(\mu,G)\in\JJ\text{ with }(\lambda,F)\not=(\mu,G)\Longrightarrow (D_{\lambda}\setminus D_{\lambda F})\cap(D_{\mu}\setminus D_{\mu G})=\emptyset.
\]
For $(\lambda,F),(\mu,G)\in\AA^n$ we write
\[
\Theta_{(\lambda,F),(\mu,G)}:=\Theta_{\XX_{D_\lambda\setminus D_{\lambda F}},\XX_{D_\mu\setminus D_{\mu G}}}\in\KK(X_n).
\]
\end{definition} 

Let $m,n\in\N^k$. To prove Proposition~\ref{comp al} we first need to show that for each $(\lambda_1,F_1),(\lambda_2,F_2)\in\AA^m$ and $(\mu_1,G_1),(\mu_2,G_2)\in\AA^n$ we have
\[
\iota_m^{m\vee n}\left(\Theta_{(\lambda_1,F_1),(\lambda_2,F_2)}\right)\iota_n^{m\vee n}\left(\Theta_{(\mu_1,G_1),(\mu_2,G_2)}\right)\in\KK(X_{m\vee n}).
\]
We do this by finding for each $(\alpha,\beta)\in\Lambda^{\min}(\lambda_2,\mu_1)$ finite subsets $\HH_{(\alpha,\beta)},\JJ_{(\alpha,\beta)}\subseteq\AA^{m\vee n}$ such that $\sqcup_{(\alpha,\beta)}\HH_{(\alpha,\beta)}$ and $\sqcup_{(\alpha,\beta)}\JJ_{(\alpha,\beta)}$ are disjoint, and
\begin{align}\label{suff eq}
\iota_m^{m\vee n}\left(\Theta_{(\lambda_1,F_1),(\lambda_2,F_2)}\right)\iota_n^{m\vee n}&\left(\Theta_{(\mu_1,G_1),(\mu_2,G_2)}\right)\\
&=\sum_{(\alpha,\beta)\in\Lambda^{\min}(\lambda_2,\mu_1)}\sum_{\substack{(\kappa,H)\in\HH_{(\alpha,\beta)} \\ (\omega,J)\in\JJ_{(\alpha,\beta)}}}\Theta_{(\kappa,H),(\omega,J)}.\nonumber
\end{align}
To find the correct $\HH_{(\alpha,\beta)}$ and $\JJ_{(\alpha,\beta)}$, we evaluate both sides of (\ref{suff eq}) on products $fg$, where $f\in C_c(\partial\Lambda^{\ge n})$ and $g\in C_c(\partial\Lambda^{\ge m\vee n-n})$. For the left-hand-side of (\ref{suff eq}) we use (\ref{char of iota map}) and Corollary~\ref{imp cor for dense subs} to factor 
\[
\iota_n^{m\vee n}\left(\Theta_{(\mu_1,G_1),(\mu_2,G_2)}\right)(fg) = \Theta_{(\mu_1,G_1),(\mu_2,G_2)}(f)g=hl,
\]
where $h\in C_c(\partial\Lambda^{\ge m})$ and $l\in C_c(\partial\Lambda^{\ge m\vee n-m})$. Then for $x\in\partial\Lambda^{\ge m\vee n}$ we have
\begin{align*}
&\iota_m^{m\vee n}\left(\Theta_{(\lambda_1,F_1),(\lambda_2,F_2)}\right)\iota_n^{m\vee n}\left(\Theta_{(\mu_1,G_1),(\mu_2,G_2)}\right)(fg)(x)\\
&= \iota_m^{m\vee n}\left(\Theta_{(\lambda_1,F_1),(\lambda_2,F_2)}\right)(hl)(x)\\
&= \Theta_{(\lambda_1,F_1),(\lambda_2,F_2)}(h)l(x)\\
&= \XX_{D_{\lambda_1}\setminus D_{\lambda_1 F_1}}(x){\langle\XX_{D_{\lambda_2}\setminus D_{\lambda_2 F_2}} ,h\rangle}_m(\sigma_m(x))l(\sigma_m(x)) \\
&= \XX_{D_{\lambda_1}\setminus D_{\lambda_1 F_1}}(x)\left(\sum_{\sigma_m(y)=\sigma_m(x)}\overline{\XX_{D_{\lambda_2}\setminus D_{\lambda_2 F_2}}(y)}h(y)\right)l(\sigma_m(x)) \\
&= 
\begin{cases}
hl(\lambda_2(0,m)\sigma_m(x)) & \text{if $x\in \left(D_{\lambda_1}\setminus D_{\lambda_1 F_1}\right)\cap\sigma_m^{-1}(D_{\lambda_2(m,d(\lambda_2))}\setminus D_{\lambda_2(m,d(\lambda_2)) F_2})$,} \\
0 & \text{otherwise.}
\end{cases} 
\end{align*}
A similar calculation to the one above gives
\begin{align*}
&hl(\lambda_2(0,m)\sigma_m(x))\\
&= \Theta_{(\mu_1,G_1),(\mu_2,G_2)}(f)g(\lambda_2(0,m)\sigma_m(x))\\
&=
\begin{cases}
fg(\mu_2(0,n)\sigma_n(\lambda_2(0,m)\sigma_m(x))) & \text{if $\lambda_2(0,m)\sigma_m(x)\in \left(D_{\mu_1}\setminus D_{\mu_1 G_1}\right)\cap$}\\
&\qquad\quad\quad \text{$\sigma_n^{-1}(D_{\mu_2(n,d(\mu_2))}\setminus D_{\mu_2(n,d(\mu_2)) G_2})$,} \\
0 & \text{otherwise.}
\end{cases} 
\end{align*}
So we label conditions
\begin{equation}\label{eq 1 for set res for comp al}
x\in \left(D_{\lambda_1}\setminus D_{\lambda_1 F_1}\right)\cap\sigma_m^{-1}(D_{\lambda_2(m,d(\lambda_2))}\setminus D_{\lambda_2(m,d(\lambda_2)) F_2}),
\end{equation}
and
\begin{equation}\label{eq 2 for set res for comp al}
\lambda_2(0,m)\sigma_m(x)\in \left(D_{\mu_1}\setminus D_{\mu_1 G_1}\right)\cap\sigma_n^{-1}(D_{\mu_2(n,d(\mu_2))}\setminus D_{\mu_2(n,d(\mu_2)) G_2}),
\end{equation}
and then we have
\begin{align}
\iota_m^{m\vee n}&\left(\Theta_{(\lambda_1,F_1),(\lambda_2,F_2)}\right)\iota_n^{m\vee n}\left(\Theta_{(\mu_1,G_1),(\mu_2,G_2)}\right)(fg)(x)\label{res for iotas}\\
&= 
\begin{cases}
fg(\mu_2(0,n)\sigma_n(\lambda_2(0,m)\sigma_m(x))) & \text{if $x$ satisfies (\ref{eq 1 for set res for comp al}) and(\ref{eq 2 for set res for comp al}),} \\
0 & \text{otherwise.}
\end{cases}\nonumber 
\end{align}
Now, for each $(\alpha,\beta)\in\Lambda^{\min}(\lambda_2,\mu_1)$, $(\kappa,H)\in\HH_{(\alpha,\beta)}$ and $(\omega,J)\in\JJ_{(\alpha,\beta)}$ we have
\begin{align*}
\Theta_{(\kappa,H),(\omega,J)}&(fg)(x)\\
&= \XX_{D_{\kappa}\setminus D_{\kappa H}}(x){\langle \XX_{D_{\omega}\setminus D_{\omega J}},fg\rangle}_{m\vee n}(\sigma_{m\vee n}(x))\\
&= \XX_{D_{\kappa}\setminus D_{\kappa H}}(x)\left(\sum_{\sigma_{m\vee n}(y)=\sigma_{m\vee n}(x)}\overline{\XX_{D_{\omega}\setminus D_{\omega J}}(y)}fg(y)\right)\\
&=
\begin{cases}
fg(\tau(0,m\vee n)\sigma_{m\vee n}(x)) & \text{if $x\in \left(D_{\kappa}\setminus D_{\kappa H}\right)\cap\sigma_{m\vee n}^{-1}(D_{\omega}\setminus D_{\omega J})$,}\\
0 & \text{otherwise.}
\end{cases}   
\end{align*}
Since $\sqcup_{(\alpha,\beta)}\HH_{(\alpha,\beta)}$ and $\sqcup_{(\alpha,\beta)}\JJ_{(\alpha,\beta)}$ are disjoint, we have
\begin{align}
&\left(\sum_{(\alpha,\beta)\in\Lambda^{\min}(\lambda_2,\mu_1)}\sum_{\substack{(\kappa,H)\in\HH \\ (\tau,J)\in\JJ}}\Theta_{(\kappa,H),(\omega,J)}\right)(fg)(x)\label{res for big sum}\\
&= 
\begin{cases}
fg(\tau(0,m\vee n)\sigma_{m\vee n}(x)) & \text{if $\displaystyle x\in\sqcup_{\substack{(\alpha,\beta) \\ (\kappa,H),(\omega,J)}}\left(D_{\kappa}\setminus D_{\rho H}\right)\cap\sigma_{m\vee n}^{-1}(D_{\omega}\setminus D_{\omega J})$,}\\
0 & \text{otherwise.}
\end{cases}\nonumber   
\end{align}
Equation~(\ref{suff eq}) now follows from (\ref{res for iotas}), (\ref{res for big sum}) and the following lemma.

\begin{lemma}\label{set res for comp al}
Let $m,n\in\N^k$, and suppose the pairs $(\lambda_1,F_1),(\lambda_2,F_2)\in\AA^m$ and $(\mu_1,G_1),(\mu_2,G_2)\in\AA^n$. Then for each pair $(\alpha,\beta)\in\Lambda^{\min}(\lambda_2,\mu_1)$ there exists finite and disjoint subsets $\HH_{(\alpha,\beta)}, \JJ_{(\alpha,\beta)}\subseteq\AA^{m\vee n}$ such that $x\in\partial\Lambda^{\ge m\vee n}$ satisfies Equations~\textup{(\ref{eq 1 for set res for comp al})} and \textup{(\ref{eq 2 for set res for comp al})} if and only if
\begin{equation}\label{single cond}
x\in\bigsqcup_{(\alpha,\beta)\in\Lambda^{\min}(\lambda_2,\mu_1)}\bigsqcup_{\substack{(\kappa,H)\in\HH_{(\alpha,\beta)} \\ (\omega,J)\in\JJ_{(\alpha,\beta)}}}\left(D_{\kappa}\setminus D_{\kappa H}\right)\cap\sigma_{m\vee n}^{-1}(D_{\omega}\setminus D_{\omega J}).
\end{equation}
Moreover, if $x$ satisfies \textup{(\ref{eq 1 for set res for comp al})} and \textup{(\ref{eq 2 for set res for comp al})} and $x\in\left(D_{\kappa}\setminus D_{\kappa H}\right)\cap\sigma_{m\vee n}^{-1}(D_{\omega}\setminus D_{\omega J})$, then we have
\[
\mu_2(0,n)\sigma_n(\lambda_2(0,m)\sigma_m(x))=\omega(0,m\vee n)\sigma_{m\vee n}(x).
\]
\end{lemma}

\begin{proof}
Recall that for $\lambda,\mu\in\Lambda$ we denote by 
\[
F(\lambda,\mu)=\{\alpha\in\Lambda:(\alpha,\beta)\in\Lambda^{\min}(\lambda,\mu)\text{ for some }\beta\in\Lambda\}.
\]
Let $(\alpha,\beta)\in\Lambda^{\min}(\lambda_2,\mu_1)$. For each $(\gamma,\delta)\in\Lambda^{\min}(\lambda_1(m,d(\lambda_1)),\lambda_2(m,d(\lambda_2))\alpha)$ we define
\begin{align*}
H_{\gamma,\alpha}:=\left(\bigcup_{\nu\in F_1}F(\lambda_1\gamma,\lambda_1\nu)\right)&\cup\left(\bigcup_{\zeta\in F_2}F(\lambda_2(m,d(\lambda_2))\alpha\delta,\lambda_2(m,d(\lambda_2))\zeta)\right)\\
&\cup\left(\bigcup_{\eta\in G_1}F(\mu_1\beta\delta,\mu_1\eta)\right),
\end{align*}
and
\[
\HH_{(\alpha,\beta)}:=\{(\lambda_1\gamma,H_{\gamma,\alpha})\in\AA^{m\vee n}:(\gamma,\delta)\in\Lambda^{\min}(\lambda_1(m,d(\lambda_1)),\lambda_2(m,d(\lambda_2))\alpha)\}.
\]
For each $(\rho,\tau)\in\Lambda^{\min}(\mu_2(n,d(\mu_2)),\mu_1(n,d(\mu_1))\beta)$ we define
\begin{align*}
J_{\rho,\beta}:=\left(\bigcup_{\xi\in G_2}F(\mu_2\rho,\mu_2\xi)\right)&\cup\left(\bigcup_{\eta\in G_1}F(\mu_1(n,d(\mu_1))\beta\tau,\mu_1(n,d(\mu_1))\eta)\right)\\
&\cup\left(\bigcup_{\zeta\in F_2}F(\lambda_2\alpha\tau,\lambda_2\zeta)\right),
\end{align*}
and
\[
\JJ_{(\alpha,\beta)}:=\{(\mu_2\rho,H_{\rho,\beta})\in\AA^{m\vee n}:(\rho,\tau)\in\Lambda^{\min}(\mu_2(n,d(\mu_2)),\mu_1(n,d(\mu_1))\beta)\}.
\]
The sets $\HH_{(\alpha,\beta)}$ and $\JJ_{(\alpha,\beta)}$ are finite sets because $\Lambda$ is finitely aligned. Since the paths in the elements of $\HH_{(\alpha,\beta)}$ are of the same length, the factorisation property ensures that each $\HH_{(\alpha,\beta)}$ is disjoint. For the same reason, each $\JJ_{(\alpha,\beta)}$ is disjoint. This explains why the second union in (\ref{single cond}) is a disjoint union. Moreover, the sets $\sqcup_{(\alpha,\beta)}\HH_{(\alpha,\beta)}$ and $\sqcup_{(\alpha,\beta)}\JJ_{(\alpha,\beta)}$ are disjoint, and hence why the first union in (\ref{single cond}) is a disjoint union.    

To prove the `only if' part of the statement, we assume $x\in\partial\Lambda^{\ge m\vee n}$ satisfies (\ref{eq 1 for set res for comp al}) and (\ref{eq 2 for set res for comp al}). We have to find pairs
\begin{align*} 
(\alpha,\beta)&\in\Lambda^{\min}(\lambda_2,\mu_1),\\
(\gamma,\delta)&\in\Lambda^{\min}(\lambda_1(m,d(\lambda_1)),\lambda_2(m,d(\lambda_2))\alpha),\text{ and}\\
(\rho,\tau)&\in\Lambda^{\min}(\mu_2(n,d(\mu_2)),\mu_1(n,d(\mu_1))\beta)
\end{align*} 
such that
\begin{itemize}
\item[(a)] $x\in D_{\lambda_1\gamma}\setminus D_{\lambda_1\gamma H_{\gamma,\alpha}}$, and
\item[(b)] $\sigma_{m\vee n}(x)\in D_{\mu_2\rho(m\vee n,d(\mu_2\rho))}\setminus D_{\mu_2\rho(m\vee n,d(\mu_2\rho)) J_{\rho,\beta}}$.
\end{itemize}
Now, we know from (\ref{eq 1 for set res for comp al}) and (\ref{eq 2 for set res for comp al}) that $\lambda_2(0,m)\sigma_m(x)\in D_{\lambda_2}\cap D_{\mu_1}$, so we take
\begin{equation}\label{the alpha beta}
(\alpha,\beta):=\left(\lambda_2(0,m)\sigma_m(x)_{\lambda_2}^{\mu_1},\lambda_2(0,m)\sigma_m(x)_{\mu_1}^{\lambda_2}\right)\in\Lambda^{\min}(\lambda_2,\mu_1).
\end{equation}
We know from (\ref{eq 1 for set res for comp al}) and (\ref{eq 2 for set res for comp al}) that $\sigma_m(x)\in D_{\lambda_1(m,d(\lambda_1))}\cap D_{\lambda_2(m,d(\lambda_2))\alpha}$, so we define $(\gamma,\delta)$ to be the pair
\begin{equation}\label{the gamma delta}
\left(\sigma_m(x)_{\lambda_1(m,d(\lambda_1))}^{\lambda_2(m,d(\lambda_2))\alpha},\sigma_m(x)_{\lambda_2(m,d(\lambda_2))\alpha}^{\lambda_1(m,d(\lambda_1))}\right)\in\Lambda^{\min}(\lambda_1(m,d(\lambda_1)),\lambda_2(m,d(\lambda_2))\alpha).
\end{equation}
We now have $\sigma_m(x)\in D_{\lambda_1(m,d(\lambda_1))\gamma}$, and this along with (\ref{eq 1 for set res for comp al}) implies that $x\in D_{\lambda_1\gamma}$. We also have
\begin{equation}\label{eq 1 for only if}
x\in D_{\lambda_1\gamma}\text{ and }x\not\in D_{\lambda_1 F_1}\Longrightarrow x\not\in D_{\lambda_1\gamma\nu'}\text{ for all }\nu'\in \bigcup_{\nu\in F_1}F(\lambda_1\gamma,\lambda_1\nu);
\end{equation}
\begin{align}
&\sigma_m(x)\in D_{\lambda_2(m,d(\lambda_2))\alpha\delta}\text{ and }\sigma_m(x)\not\in D_{\lambda_2(m,d(\lambda_2)) F_2}\nonumber\\
& \Longrightarrow \sigma_m(x)\not\in D_{\lambda_2(m,d(\lambda_2))\alpha\delta\zeta'}\text{ for all }\zeta'\in\bigcup_{\zeta\in F_2}F(\lambda_2(m,d(\lambda_2))\alpha\delta,\lambda_2(m,d(\lambda_2))\zeta)\nonumber\\
&\Longleftrightarrow \sigma_m(x)\not\in D_{\lambda_1(m,d(\lambda_1))\gamma\zeta'}\text{ for all }\zeta'\in\bigcup_{\zeta\in F_2}F(\lambda_2(m,d(\lambda_2))\alpha\delta,\lambda_2(m,d(\lambda_2))\zeta)\nonumber\\
&\Longleftrightarrow x\not\in D_{\lambda_1\gamma\zeta'}\text{ for all }\bigcup_{\zeta\in F_2}F(\lambda_2(m,d(\lambda_2))\alpha\delta,\lambda_2(m,d(\lambda_2))\zeta);\label{eq 2 for only if}
\end{align}
and
\begin{align}
\lambda_2(0,m)\sigma_m(x)\in & D_{\lambda_2\alpha\delta}\text{ and }\lambda_2(0,m)\sigma_m(x)\not\in D_{\mu_1G_1}\nonumber\\
&\Longrightarrow \lambda_2(0,m)\sigma_m(x)\not\in D_{\lambda_2\alpha\delta\eta'}\text{ for all }\eta'\in \bigcup_{\eta\in G_1}F(\mu_1\beta\delta,\mu_1\eta)\nonumber\\
&\Longleftrightarrow \sigma_m(x)\not\in D_{\lambda_2(m,d(\lambda_2))\alpha\delta\eta'}\text{ for all }\eta'\in \bigcup_{\eta\in G_1}F(\mu_1\beta\delta,\mu_1\eta)\nonumber\\
&\Longleftrightarrow \sigma_m(x)\not\in D_{\lambda_1(m,d(\lambda_1))\gamma\eta'}\text{ for all }\eta'\in \bigcup_{\eta\in G_1}F(\mu_1\beta\delta,\mu_1\eta)\nonumber\\
&\Longleftrightarrow x\not\in D_{\lambda_1\gamma\eta'}\text{ for all }\eta'\in \bigcup_{\eta\in G_1}F(\mu_1\beta\delta,\mu_1\eta).\label{eq 3 for only if}
\end{align}
It follows from (\ref{eq 1 for only if}), (\ref{eq 2 for only if}) and (\ref{eq 3 for only if}) that $x\not\in D_{\lambda_1\gamma H_{\gamma,\alpha}}$, and so (a) is satisfied.

We have $\sigma_n(\lambda_2(0,m)\sigma_m(x))\in D_{\mu_1(n,d(\mu_1))\beta}$, and it follows from (\ref{eq 2 for set res for comp al}) that $\sigma_n(\lambda_2(0,m)\sigma_m(x))\in D_{\mu_2(n,d(\mu_2))}$. So we take
\begin{align}
(\rho,\tau)&:= \left(\sigma_n(\lambda_2(0,m)\sigma_m(x))_{\mu_2(n,d(\mu_2))}^{\mu_1(n,d(\mu_1))\beta},\sigma_n(\lambda_2(0,m)\sigma_m(x))_{\mu_1(n,d(\mu_1))\beta}^{\mu_2(n,d(\mu_2))}\right)\label{the rho tau}\\
&\in\Lambda^{\min}(\mu_2(n,d(\mu_2)),\mu_1(n,d(\mu_1))\beta),\nonumber
\end{align}
and we have
\begin{align*}
\sigma_n(\lambda_2(0,m)\sigma_m(x))\in&  D_{\mu_2(n,d(\mu_2))\rho}\\
& \Longrightarrow \sigma_{m\vee n}(x)=\sigma_{m\vee n}(\lambda_2(0,m)\sigma_m(x))\in D_{\mu_2\rho(m\vee n,d(\mu_2\rho))}.
\end{align*}
Suppose for contradiction that there exists $\xi\in G_2$ and a pair $(\xi',\xi'')$ in the set $\Lambda^{\min}(\mu_2\rho,\mu_2\xi)$ with $\sigma_{m\vee n}(x)\in D_{\mu_2\rho(m\vee n,d(\mu_2\rho))\xi'}$. Then it follows from (\ref{the rho tau}) that
\begin{align*}
\sigma_n(\lambda_2(0,m)\sigma_m(x))&=\sigma_n(\lambda_2(0,m)\sigma_m(x))(0,m\vee n -n)\sigma_{m\vee n}(x)\\
&=\mu_2(n,d(\mu_2))\rho(0,m\vee n -n)\sigma_{m\vee n}(x)\\
&= \mu_2\rho(n,m\vee n)\sigma_{m\vee n}(x)\\
&\in D_{\mu_2(n,d(\mu_2))\rho\xi'}\\
&=D_{\mu_2(n,d(\mu_2))\xi\xi''}\\
&\subseteq D_{\mu_2(n,d(\mu_2)) G_2}.
\end{align*}
This contradicts Equation~(\ref{eq 2 for set res for comp al}), and so we must have
\begin{equation}\label{eq 1 for if}
\sigma_{m\vee n}(x)\not\in D_{\mu_2\rho(m\vee n,d(\mu_2\rho))\xi'}\text{ for all }\xi'\in \bigcup_{\xi\in G_2}F(\mu_2\rho,\mu_2\xi).
\end{equation}
Similar arguments show that
\begin{equation}\label{eq 2 for if}
\sigma_{m\vee n}(x)\not\in D_{\mu_2\rho(m\vee n,d(\mu_2\rho))\eta'},\text{ for all $\eta'\in \bigcup_{\eta\in G_1}F(\mu_1(n,d(\mu_1))\beta\tau,\mu_1(n,d(\mu_1))\eta)$},
\end{equation}
and
\begin{equation}\label{eq 3 for if}
\sigma_{m\vee n}(x)\not\in D_{\mu_2\rho(m\vee n,d(\mu_2\rho))\eta'},\text{ for all $\zeta'\in \bigcup_{\zeta\in F_2}F(\lambda_2\alpha\tau,\lambda_2\zeta)$}.
\end{equation}
It follows from (\ref{eq 1 for if}), (\ref{eq 2 for if}) and (\ref{eq 3 for if}) that $\sigma_{m\vee n}(x)\not\in D_{\mu_2\rho(m\vee n,d(\mu_2\rho)) J_{\rho,\beta}}$, and so (b) is satisfied.

To prove the `if' part of the statement, we assume there exists
\begin{align*} 
(\alpha,\beta)&\in\Lambda^{\min}(\lambda_2,\mu_1),\\
(\gamma,\delta)&\in\Lambda^{\min}(\lambda_1(m,d(\lambda_1)),\lambda_2(m,d(\lambda_2))\alpha),\text{ and}\\
(\rho,\tau)&\in\Lambda^{\min}(\mu_2(n,d(\mu_2)),\mu_1(n,d(\mu_1))\beta),
\end{align*}
such that
\[
x\in \left(D_{\lambda_1\gamma}\setminus D_{\lambda_1\gamma H_{\gamma,\alpha}}\right)\cap\sigma_{m\vee n}^{-1}(D_{\mu_2\rho(m\vee n, d(\mu_2\rho))}\setminus D_{\mu_2\rho(m\vee n, d(\mu_2\rho)) J_{\rho,\beta}}).
\]
We have
\[
x\in D_{\lambda_1\gamma}\setminus D_{\lambda_1\gamma H_{\gamma,\alpha}}
\Longrightarrow  x\in D_{\lambda_1}\setminus D_{\lambda_1 F_1},
\]      
and
\begin{align*}
x\in  D_{\lambda_1\gamma}\setminus D_{\lambda_1\gamma H_{\gamma,\alpha}}
&\Longrightarrow \sigma_m(x)\in D_{\lambda_1(m,d(\lambda_1))\gamma}\setminus D_{\lambda_1(m,d(\lambda_1))\gamma H_{\gamma,\alpha}}\\
&\Longleftrightarrow \sigma_m(x)\in D_{\lambda_2(m,d(\lambda_2))\alpha\delta}\setminus D_{\lambda_2(m,d(\lambda_2))\alpha\delta H_{\gamma,\alpha}}\\
&\Longleftrightarrow \sigma_m(x)\in D_{\lambda_2(m,d(\lambda_2))}\setminus D_{\lambda_2(m,d(\lambda_2)) F_2}.
\end{align*}
So (\ref{eq 1 for set res for comp al}) is satisfied. We have
\begin{align*}
x\in D_{\lambda_1\gamma}\setminus D_{\lambda_1\gamma H_{\gamma,\alpha}} &\Longrightarrow \sigma_m(x)\in D_{\lambda_1(m,d(\lambda_1))\gamma}\setminus D_{\lambda_1(m,d(\lambda_1))\gamma H_{\gamma,\alpha}}\\
&\Longleftrightarrow \sigma_m(x)\in D_{\lambda_2(m,d(\lambda_2))\alpha\delta}\setminus D_{\lambda_2(m,d(\lambda_2))\alpha\delta H_{\gamma,\alpha}}\\
&\Longrightarrow \lambda_2(0,m)\sigma_m(x)\in D_{\lambda_2\alpha\delta}\setminus D_{\lambda_2\alpha\delta H_{\gamma,\alpha}}\\
&\Longleftrightarrow \lambda_2(0,m)\sigma_m(x)\in D_{\mu_1\beta\delta}\setminus D_{\mu_1\beta\delta H_{\gamma,\alpha}}\\
&\Longleftrightarrow \lambda_2(0,m)\sigma_m(x)\in D_{\mu_1}\setminus D_{\mu_1 G_1}.
\end{align*}
We have
\begin{align*}
x\in D_{\lambda_1\gamma}\Longrightarrow \lambda_2(0,m)\sigma_m(x)(n,m\vee n) &= \big(\lambda_2(0,m)\lambda_1(m,d(\lambda_1))\gamma\big)(n,m\vee n)\\
&= \big(\lambda_2(0,m)\lambda_2(m,d(\lambda_2))\alpha\delta\big)(n,m\vee n)\\
&= \lambda_2\alpha\delta(n,m\vee n)\\
&= \lambda_2\alpha(n,m\vee n)\\
&= \mu_1\beta(n,m\vee n)\\
&= \big(\mu_1(n,d(\mu_1))\beta\big)(n,m\vee n)\\
&= \big(\mu_1(n,d(\mu_1))\beta\tau\big)(n,m\vee n)\\
&= \big(\mu_2(n,d(\mu_2))\rho\big)(n,m\vee n).
\end{align*}
It follows that
\begin{align*}
\sigma_n(\lambda_2(0,m)\sigma_m(x)) &= \big(\lambda_2(0,m)\sigma_m(x)\big)(n,m\vee n)\sigma_{m\vee n}(\lambda_2(0,m)\sigma_m(x))\\
&= \big(\lambda_2(0,m)\sigma_m(x)\big)(n,m\vee n)\sigma_{m\vee n}(x)\\
&= \big(\mu_2(n,d(\mu_2))\rho\big)(n,m\vee n)\sigma_{m\vee n}(x),
\end{align*}
and then we have
\begin{align*}
\sigma_{m\vee n}(x)\in & D_{\mu_2\rho(m\vee n,d(\mu_2\rho))}\setminus D_{\mu_2\rho(m\vee n,d(\mu_2\rho)) J_{\rho,\beta}}\\
&\Longrightarrow \sigma_n(\lambda_2(0,m)\sigma_m(x))\in D_{\mu_2\rho(n,d(\mu_2\rho))}\setminus D_{\mu_2\rho(n,d(\mu_2\rho)) J_{\rho,\beta}}\\ 
&\Longrightarrow \sigma_n(\lambda_2(0,m)\sigma_m(x))\in D_{\mu_2(n,d(\mu_2))}\setminus D_{\mu_2(n,d(\mu_2)) G_2}.
\end{align*} 
So (\ref{eq 2 for set res for comp al}) is satisfied.

To prove the final part of the result, recall that, given $x\in\partial\Lambda^{\ge m\vee n}$ satisfying (\ref{eq 1 for set res for comp al}) and (\ref{eq 2 for set res for comp al}), we have the following formula for the pair $(\rho,\tau)$ in the set $\Lambda^{\min}(\mu_2(n,d(\mu_2)),\mu_1(n,d(\mu_1))\beta)$:
\[
(\rho,\tau)=\left(\sigma_n(\lambda_2(0,m)\sigma_m(x))_{\mu_2(n,d(\mu_2))}^{\mu_1(n,d(\mu_1))\beta},\sigma_n(\lambda_2(0,m)\sigma_m(x))_{\mu_1(n,d(\mu_1))\beta}^{\mu_2(n,d(\mu_2))}\right).
\]
We then have
\begin{align*}
\mu_2(0,n)\sigma_n(\lambda_2(0,m)\sigma_m(x))&=\mu_2(0,n)\big(\sigma_n(\lambda_2(0,m)\sigma_m(x))\big)(0,m\vee n-n)\sigma_{m\vee n}(x)\\
&= \mu_2(0,n)\big(\mu_2(n,d(\mu_2))\rho\big)(0,m\vee n-n)\sigma_{m\vee n}(x)\\
&= \mu_2\rho(0,m\vee n)\sigma_{m\vee n}(x).\qedhere
\end{align*}
\end{proof}

\begin{proof}[Proof of Proposition~\ref{comp al}]
We have already established Equation~(\ref{suff eq}). Since $\Lambda$ is finitely aligned, the sums in (\ref{suff eq}) are finite, and so 
\[
\iota_m^{m\vee n}\left(\Theta_{(\lambda_1,F_1),(\lambda_2,F_2)}\right)\iota_n^{m\vee n}\left(\Theta_{(\mu_1,G_1),(\mu_2,G_2)}\right)\in\KK(X_{m\vee n}),
\]
for every $m,n\in\N^k$, $(\lambda_1,F_1),(\lambda_2,F_2)\in\AA^m$ and $(\mu_1,G_1),(\mu_2,G_2)\in\AA^n$. It then follows from Proposition~\ref{res for X m} that $\iota_m^{m\vee n}\left(\Theta_{x_1,x_2}\right)\iota_n^{m\vee n}\left(\Theta_{y_1,y_2}\right)\in\KK(X_{m\vee n})$, for every $x_1,x_2\in X_m$ and $y_1,y_2\in X_n$. Hence, $\iota_m^{m\vee n}\left(S\right)\iota_n^{m\vee n}\left(T\right)\in\KK(X_{m\vee n})$, for every $S\in\KK(X_m)$ and $(T\in\KK(X_n)$.
\end{proof}

\begin{acknowledgements} 
The author thanks Iain Raeburn for the many helpful discussions while 
this research was being conducted. This research was supported by the Australian Research Council
\end{acknowledgements}


\end{document}